%% file: connex.tex
\documentclass[11pt]{article}
\usepackage[latin1]{inputenc}
\usepackage{amsmath,amsfonts,amscd,amssymb}
\usepackage[all]{xy}
\usepackage{graphicx,epsfig}
\usepackage{color}
\usepackage[]{hyperref}
\usepackage{makeidx}
\newtheorem{thm}{Theorem}[section]
\newtheorem{cor}{Corollary}[section]
\newtheorem{lem}{Lemma}[section]
\newtheorem{prop}{Proposition}[section]
\newtheorem{deft}{Definition}[section]
\newtheorem{rek}{Remark}[section]
\newtheorem{conj}{Conjecture}
\let\noi=\noindent
\let\ss=\underline
\let\sur=\overline

\newcommand{\dsp}{\displaystyle}
\newcommand{\void}{\emptyset}
\newcommand{\rr}{\mbox{rat.rk}}
\def\N{\mathbb{N}} 
\def\Z{\mathbb{Z}}  
\def\R{\mathbb{R}}  
\def\Q{\mathbb{Q}}
\def\F{\mathbb{F}}
\def\sper{\mbox{Sper}}

\def\notin{\mbox{$\in$ \hspace{-.8em}/}} 
\def\sgn{\mbox{sgn}} 

\title{ON CONNECTEDNESS OF SETS IN THE REAL SPECTRA OF POLYNOMIAL RINGS}
\author{F. Lucas\\
D\'{e}partement de Math\'{e}matiques\\
Universit\'{e} d'Angers\\
2, bd Lavoisier\\
49045 Angers Cedex, France \and
J. J. Madden\\
Department of Mathematics\\
Louisiana State University at Baton Rouge,\\
Baton Rouge, LA, USA \and
D. Schaub\\
D\'{e}partement de Math\'{e}matiques\\
Universit\'{e} d'Angers\\
2, bd Lavoisier\\
49045 Angers Cedex, France \and
M. Spivakovsky\\
Insitut de Math\'{e}matiques de Toulouse\\
Universit\'e Paul Sabatier\\
118, route de Narbonne\\
31062 Toulouse Cedex 9, France.}

\date{}

\begin{document}

\maketitle

\begin{abstract} Let $R$ be a real closed field.
The Pierce-Birkhoff conjecture says that any piecewise polynomial
function $f$ on $R^n$ can be obtained from the polynomial ring
$R[x_1,\dots,x_n]$ by iterating the operations of maximum and
minimum. The purpose of this paper is threefold. First, we state a new
conjecture, called the Connectedness conjecture, which asserts, for every pair of points $\alpha,\beta\in\sper\
R[x_1,\dots,x_n]$, the existence of connected sets in the real spectrum of
$R[x_1,\dots,x_n]$, satisfying certain conditions. We prove that the
Connectedness conjecture implies the Pierce-Birkhoff conjecture.

Secondly, we construct a class of connected sets in the real spectrum
which, though not in itself enough for the proof of the
Pierce-Birkhoff conjecture, is the first and simplest example of the
sort of connected sets we really need, and which constitutes the first
step in our program for a proof of the Pierce-Birkhoff conjecture in
dimension greater than 2.

Thirdly, we apply these ideas to give two proofs of the Connectedness conjecture (and hence also of the Pierce--Birkhoff
conjecture in the abstract formulation) in the special case when one of the two points $\alpha,\beta\in\sper\
R[x_1,\dots,x_n]$ is monomial. One of the proofs is elementary while the other consists in deducing the (monomial)
Connectedness conjecture as an immediate corollary of the main connectedness theorem of this paper.

\end{abstract}

\section{Introduction}
\label{In}

All the rings in this paper will be commutative with 1. Throughout this paper, $R$ will denote a real closed field and $A$
the polynomial ring $R[x_1,\dots,x_n]$, unless otherwise specified.

The Pierce-Birkhoff conjecture asserts that any piecewise-polynomial
function $f:R^n\rightarrow R$ can be expressed as a maximum of minima
of a finite family of polynomials (see below for the definitions and a
precise statement of the conjecture). This paper is the first step in our program
for a proof of the Pierce-Birkhoff
conjecture in its full generality (the best results up to now are due to Louis
Mah\'e \cite{Mah}, who proved the conjecture for $n=2$, as well as some partial results for $n=3$).

We start by stating the Pierce--Birkhoff conjecture in its original
form as it was first stated by M. Henriksen and J. Isbell in the early nineteen sixties.
\begin{deft}\label{pw} A function $f:R^n\to R$ is said to be \textbf{piecewise
  polynomial} if $R^n$ can be covered by a finite collection of closed
  semi-algebraic sets $P_i$ such that for each $i$ there exists a
  polynomial $f_i\in A$ satisfying
  $\left.f\right|_{P_i}=\left.f_i\right|_{P_i}$. Given a piecewise polynomial function $f$, we say
  that $f$ is \textbf{defined by} $r$ polynomials if $r$ is the number of distinct polynomials among the $f_i$ above.
\end{deft}
Clearly, any piecewise polynomial function is continuous. Piecewise polynomial functions
form a ring, containing $A$, which is denoted by $PW(A)$.\medskip

On the other hand, one can consider the (lattice-ordered) ring of all
the functions obtained from $A$ by iterating the operations of $\sup$
and $\inf$. Since applying the operations of sup and inf to
polynomials produces functions which are piecewise polynomial, this
ring is contained in $PW(A)$ (the latter ring is closed under $\sup$
and $\inf$). It is natural to ask whether the two rings coincide. The
precise statement of the conjecture is:
\begin{conj}\textnormal{\textbf{(Pierce-Birkhoff)}}\label{PB} If
  $f:R^n\to R$ is in $PW(A)$, then there exists a finite family of
  polynomials $g_{ij}\in A$ such that
  $f=\sup\limits_i\inf\limits_j(g_{ij})$ (in other words, for all
  $x\in R^n$, $f(x)=\sup\limits_i\inf\limits_j(g_{ij}(x))$).
\end{conj}

A few words about the history of the problem. A question suggestive of the present Pierce--Birkhoff conjecture
appeared in the 1956 paper of Birkhoff and Pierce \cite{BiPi}, where f-rings were
first defined. The question seems to have been phrased somewhat
carelessly, and the intended meaning of the text is not easy to
understand.

According to Mel Henriksen, he and Isbell attempted to prove the
conjecture as it is now formulated while working on f-rings in the
early nineteen sixties \cite{HenIsb}. Since they obtained no
significant results about it,
they never mentioned it explicitly in print.

In the early nineteen eighties, when real algebraic geometry was
becoming active, Isbell described the question to several mathematicians at several
meetings. Efroymson and Mah\'e took up the challenge and worked on it.

Isbell did not hesitate to read the 1956 question posed by Birkhoff
and Pierce quite generously. He seems to be the one responsible for
attaching the name ``Pierce--Birkhoff conjecture'' to the problem. He
believed that Birkhoff and Pierce intended to ask the
very question we now call the Pierce--Birkhoff conjecture.

It is not clear how or why the names Pierce and Birkhoff came to
appear in reverse alphabetical order, rather than the alphabetical
order that was used in their 1956 paper. In light of the history, it
might be more accurate to call the question the
``Pierce-Birkhoff-Isbell-Henriksen conjecture''.

In 1989 J.J. Madden \cite{Mad1} reformulated this conjecture in terms
of the real spectrum of $A$ and separating ideals. We now recall
Madden's formulation together with the relevant definitions.

Let $B$ be a ring. A point $\alpha$ in the real spectrum of $B$ is, by
definition, the data of a prime ideal $\mathfrak{p}$ of $B$, and a
total ordering $\leq$ of the quotient ring $B/\mathfrak{p}$, or,
equivalently, of the field of fractions of $B/\mathfrak{p}$.  Another
way of defining the point $\alpha$ is as a homomorphism from $B$ to a
real closed field, where two homomorphisms are identified if they have
the same kernel $\mathfrak{p}$ and induce the same total ordering on
$B/\mathfrak{p}$.

The ideal $\mathfrak{p}$ is called the support of $\alpha$ and denoted
by $\mathfrak{p}_\alpha$, the quotient ring $B/\mathfrak{p}_{\alpha}$
by $B[\alpha]$, its field of fractions by $B(\alpha)$ and the real
closure of $B(\alpha)$ by $k(\alpha)$. The total ordering of
$B(\alpha)$ is denoted by $\leq_\alpha$. Sometimes we write
$\alpha=(\mathfrak{p}_\alpha,\leq_\alpha)$.
\begin{deft} The real spectrum of $B$, denoted by Sper $B$, is the
  collection of all pairs $\alpha=(\mathfrak{p}_\alpha,\leq_\alpha)$,
  where $\mathfrak{p}_\alpha$ is a prime ideal of $B$ and
  $\leq_\alpha$ is a total ordering of $B/\mathfrak{p}_{\alpha}$.
\end{deft}
The real spectrum $\sper\ B$ is endowed with two natural topologies.
The first one, called the {\bf spectral (or Harrison) topology}, has basic open sets of the form
$$
U(f_1,\ldots,f_k)=\{\alpha\ |\ f_1(\alpha)>0,\ldots,f_k(\alpha)>0\}
$$
with $f_1,...,f_k\in B$. Here and below, we commit the following
standard abuse of notation: for an element $f\in B$, $f(\alpha)$
stands for the natural image of $f$ in $B[\alpha]$ and the inequality
$f(\alpha)>0$ really means $f(\alpha)>_\alpha0$.

The second is the \textbf{constructible topology} whose basic open sets
are of the form
$$
V(f_1,\ldots,f_k,g)=\{\alpha\ |\
f_1(\alpha)>0,\ldots,f_k(\alpha)>0,g(\alpha)=0\},
$$
where $f_1,...,f_n,g\in B$. Boolean combinations of sets of the form $V(f_1,\ldots,f_n,g)$ are
called \textbf{constructible sets} of $\sper\ B$.

\begin{prop} (\cite{BCR}, Proposition 7.1.12, p.114) Let $B$ be any ring, then $\sper(B)$ is compact for the constructible
topology.
\end{prop}

\begin{rek} Since the spectral topology is coarser than the constructible topology, $\sper(B)$ is also compact for the
spectral topology.
\end{rek}

Denote by Maxr$(A)$ the set of points $\alpha \in \sper(A)$ such that
$\mathfrak{p}_\alpha$ is a maximal ideal of $A$. We view Maxr$(A)$ as
a topological subspace of $\sper(A)$ with the spectral (respectively,
constructible) topology. We may naturally identify $R^n$
with Maxr$(A)$: a point $(a_1,\ldots,a_n)\in R^n$ corresponds to the point
$\alpha=(\mathfrak{p}_\alpha,\leq_\alpha)\in\sper(A)$, where $\mathfrak{p}_\alpha$
is the maximal ideal
$$
\mathfrak{p}_\alpha=(x_1-a_1,\ldots,x_n-a_n)
$$
and $\leq_\alpha$ is the unique order on $R$. The spectral topology on $\sper(A)$ induces
the euclidean topology on $R^n$, while the constructible topology
induces the topology of $R^n$ whose base consists of the
semi-algebraic sets. The injection $R^n=\mbox{Maxr}(A)\hookrightarrow\sper(A)$ induces a one-to-one correspondence between
the semi-algebraic sets of $R^n$ and the constructible sets of $\sper(A)$: the set
$S(f_1,\ldots,f_k,g)=\{(a_1,\ldots,a_n)\in R^n \ |\
f_1(a_1,\ldots,a_n) > 0,\ldots,f_k(a_1,\ldots,a_n)>0,g(a_1,\ldots,a_n)=0 \}$ corresponds to the set
$$
V(f_1,\ldots,f_k,g)=\{\alpha\ |\ f_1(\alpha)>0,\ldots,f_k(\alpha)>0,g(\alpha)=0\}.
$$
For a general semi-algebraic set $S$ in $R^n$, the corresponding constructible set
of $\sper(A)$ will be denoted by $ \tilde{S}$. Conversely, we
sometimes start with a constructible set $\tilde{S}$ in $\sper(A)$ and
denote by $S$ its semi-algebraic counterpart in $R^n$.
\medskip

\noi\textbf{Notation:} We will use the notation $\coprod$ to denote the
disjoint union of sets.
\medskip

Recall Definition 2.4.2 and Proposition 7.5.1 of \cite{BCR} (p. 31 and p. 130, respectively):
\begin {deft} Let $S$ be a semi-algebraic subset of $R^n$. The set $S$
  is called semi-algebraically connected if, for any decomposition
  $S= S_1 \coprod S_2$, with $S_1,S_2$ semi-algebraic and relatively
  closed in $S$ for the euclidean topology, we have either $S_1 = S$
  or $S_2=S$.
\end{deft}

\begin{prop}\label{tilde} A semi-algebraic set $S$ in $R^n$ is semi-algebraically
connected if and only if $\tilde{S}$ is connected in the spectral topology.
If $S_1,\ldots,S_k$ are the semi-algebraically connected components of $S$,
then $\tilde{S}_1,\ldots,\tilde{S}_k$ are the connected components of
$\tilde{S}$.
\end{prop}

\noi\textbf{Warning:} Semi-algebraic connectedness is the same as the usual
connectedness if $R=\R$, but not in general. For a general real closed
field $R$, usual connectedness in the euclidean topology is a very
restrictive notion and will not be used in this paper.

\begin{prop} \label{qualite}(\cite{BCR}, Theorems 2.4.4 and 7.5.1). Let $A=R[x_1,\dots,x_n]$ and let $X$ be a
  constructible subset of $\sper(A)$. Then $X$ can be written as
\begin{equation}
  X =  \coprod\limits_{i=1}^r X_i\label{eq:decomp2}
\end{equation}
 where each $X_i$ is constructible and connected
  in the spectral topology. In other words, $X$ has finitely many
  connected components in the spectral topology, each of which is itself constructible.
\end{prop}
\begin{rek}\label{fingen} Let $B$ be a finitely generated $R$-algebra. Then Proposition \ref{qualite} holds also for $\sper\ A$ replaced
with $\sper\ B$. Indeed, $B$ is a homomorphic image of a polynomial ring of the form $R[x_1,\dots,x_n]$, so $\sper\ B$ is a
closed subspace of $\sper\ R[x_1,\dots,x_n]$. Then a constructible subset $X$ of $\sper\ B$ is also constructible viewed
as a subset of $\sper\ R[x_1,\dots,x_n]$ and the decomposition (\ref{eq:decomp2}) has the properties required in the
Proposition regardless of whether one views the $X_i$ as subsets of $\sper\ B$ or of $\sper\ R[x_1,\dots,x_n]$.
\end{rek}
Next, we recall the notion of \textbf{separating ideal}, introduced by
Madden in \cite{Mad1}.
\begin{deft} Let $B$ be a ring. For $\alpha,\beta\in\mbox{Sper}\ B$,
  the \textbf{separating ideal} of $\alpha$ and $\beta$, denoted by
  $<\alpha,\beta>$, is the ideal of $B$ generated by all the elements
  $f\in B$ which change sign between $\alpha$ and $\beta$, that is,
  all the $f$ such that $f(\alpha)\geq0$ and $f(\beta)\leq0$.
\end{deft}
We will need the following basic properties of the separating ideal, proved in \cite{Mad1}:
\begin{prop}\label{sepideal} Let the notation be as above. We
have:

(1) $<\alpha,\beta>$ is both a $\nu_\alpha$-ideal and a $\nu_\beta$-ideal.

(2) $\alpha$ and $\beta$ induce the same ordering on $\frac B{<\alpha,\beta>}$ (in particular, the set of
$\nu_\alpha$-ideals containing $<\alpha,\beta>$ coincides with the set of $\nu_\beta$-ideals containing $<\alpha,\beta>$).

(3) $<\alpha,\beta>$ is the smallest ideal (in the sense of inclusion), satisfying (1) and (2).

(4) If $\alpha$ and $\beta$ have no common specialization then $<\alpha,\beta>=B$.
\end{prop}
Let $f$ be a piecewise polynomial function on $R^n$ and take a point
$\alpha\in\sper\ A$. Let the notation be as in Definition \ref{pw}.
The covering $R^n=\bigcup\limits_iP_i$ induces a corresponding
covering $\sper\ A=\bigcup\limits_i\tilde P_i$ of the real spectrum.
Pick and fix an $i$ such that $\alpha\in\tilde P_i$. We set
$f_\alpha:=f_i$. We refer to $f_\alpha$ as {\bf a local polynomial
representative of $f$ at $\alpha$}. In general, the choice of $i$ is
not uniquely determined by $\alpha$.  Implicit in the notation
$f_\alpha$ is the fact that one such choice has been made.

In \cite{Mad1}, Madden reduced the Pierce--Birkhoff conjecture to a
purely local statement about separating ideals and the real spectrum.
Namely, he showed that the Pierce-Birkhoff conjecture is equivalent to
\begin{conj}\textnormal{\textbf{(Pierce-Birkhoff conjecture, the
abstract version)}}\label{PBS} Let $f$ be a piecewise polynomial
function and $\alpha,\beta$ points in $\sper\ A$. Let $f_\alpha\in
A$ be a local representative of $f$ at $\alpha$ and $f_\beta\in A$ a
local representative of $f$ at $\beta$. Then
$f_\alpha-f_\beta\in<\alpha,\beta>$.
\end{conj}
The following slightly weaker conjecture has proved to be an extremely useful stepping stone on the
way to Pierce--Birkhoff (we give the contrapositive of Madden's original
statement, since it is better adapted to our needs).
\begin{conj}\textnormal{\textbf{(the Separation
conjecture)}}\label{Sep} Let $g\in A$ and
let $\alpha,\beta\in\sper\ A$ be two points such that
$g\notin<\alpha,\beta>$. Then $\alpha$ and $\beta$ lie in the same
connected component of the set $\sper\ A\setminus\{g=0\}$.
\end{conj}
\begin{prop} The Separation conjecture is equivalent to the special case of the Pierce--Birkhoff conjecture in which the
piecewise-polynomial function $f$ is defined by two polynomials.
\end{prop}
\noi{\it Proof:} First, assume the Pierce--Birkhoff conjecture is true for any piecewise-polynomial $f$ defined by two
polynomials. Let $g$ be as in the Separation conjecture. We argue by contradiction. Assume that $\alpha$ and $\beta$ lie in
two different connected components of $\sper\ A\setminus\{g=0\}$. Let $f$ be the piecewise polynomial function which is
equal to $g$ on the connected component of $\sper\ A\setminus\{g=0\}$ containing $\alpha$
and $f=0$ elsewhere. Then by the Pierce--Birkhoff conjecture we have $g-0=g\in<\alpha,\beta>$, which gives the desired
contradiction.

Conversely, assume that the Separation conjecture holds. Let $f$ be a piecewise polynomial function defined by two
distinct polynomials, $f_1$ and $f_2$. Let $P_1$ and $P_2$ be closed semialgebraic sets such that
$\left.f\right|_{P_i}=\left.f_i\right|_{P_i}$, $i\in\{1,2\}$. Let $g=f_1-f_2$. Then the interior of $P_1$ is disjoint from
the interior of $P_2$ and each is a union of connected components of $\sper\ A\setminus\{g=0\}$. Pick a point $\alpha$ in
the interior of $P_1$ and $\beta$ in the interior of $P_2$. Then $\alpha$ and $\beta$ lie in two different connected
components of $\sper\ A\setminus\{g=0\}$, hence $f_1-f_2\in<\alpha,\beta>$ by the separation conjecture. Thus the
Pierce--Birkhoff conjecture holds for $f$. \hfill$\Box$\medskip

We now state
\begin{conj}\textnormal{\textbf{(the Connectedness
conjecture)}}\label{conn} Let $\alpha,\beta\in\sper\ A$ and let
  $g_1,\dots,g_s$ be a finite collection of elements of $A\setminus<\alpha,\beta>$. Then there exists a connected set
  $C\subset\sper\ A$ such that $\alpha,\beta\in C$ and
  $C\cap\{g_i=0\}=\emptyset$ for $i\in\{1,\dots,s\}$ (in other words,
  $\alpha$ and $\beta$ belong to the same connected component of the
  set $\sper\ A\setminus \{g_1\dots g_s=0\}$).
\end{conj}
One advantage of the Connectedness conjecture is that it is a statement about polynomials which makes no mention of
piecewise polynomial functions.

The purpose of this paper is threefold. First, we prove
(\S\ref{implies}) that the Connectedness conjecture implies the
Pierce-Birkhoff conjecture. This reduces the Pierce-Birkhoff
conjecture to constructing, for each $\alpha,\beta\in\sper\ A$,
a connected set in $\sper\ A$ having certain properties.

Secondly, we construct a class of connected subsets of $\sper\ A$
which, though not in itself enough for the proof of the
Pierce-Birkhoff conjecture, is the first and simplest example of the
sort of connected sets we really need, and which constitutes the first
step in our program for a proof of the Pierce-Birkhoff conjecture in
dimension greater than $2$. The precise relation of the main
connectedness theorem of the present paper (Theorem \ref{con2}) to the
Pierce-Birkhoff conjecture, that is, the part of the program which is
left for the subsequent papers, is explained in more detail later in
this Introduction. In the forthcoming paper \cite{LMSS}, the next
paper of the series, we plan to prove the separation conjecture for $n=3$ by
constructing higher level connected sets in the real spectrum, based
on the prototype constructed in the present paper.

Thirdly, we give two proofs of the Connectedness conjecture, one elementary
and one invoking the main connectedness theorem of this paper, in the special case when one of $\alpha$ and $\beta$ is
monomial (this is explained in more detail later in this Introduction).
\medskip

The rest of this paper is organized as follows. In \S\ref{Cu} we
define the valuation $\nu_\alpha$ of $B(\alpha)$ associated to a point
$\alpha$ of the real spectrum. The valuation $\nu_\alpha$ is defined
using the order $\leq_\alpha$. It has following properties:

(1) $\nu_\alpha(B[\alpha]) \geq 0$

(2) If $B$ is an $R$-algebra then for any positive elements $y,z \in B(\alpha)$,
\begin{equation}
\nu_\alpha(y) < \nu_\alpha(z) \implies y > Nz,\ \forall N \in
R. \label{val2}
\end{equation}

\noi We also explain the geometric interpretation of a point of the
real spectrum as a semi-curvette.\medskip

\noi\textbf{Notation:} To simplify the notation, we will write $x$ in place
of $(x_1,\ldots,x_n)$.\medskip

For a point $\alpha$ of the real spectrum, let $\Gamma_\alpha$ denote the value group of $\nu_\alpha$. We say that $\alpha$
is \textbf{monomial} if it is given by a semi-curvette of the form $t\rightarrow(c_1t^{a_1},\dots,c_nt^{a_n})$, where
$c_j\in R$ and the $a_j$ are strictly positive elements, $\Q$-linearly independent elements of $\Gamma_\alpha$. Note that
if $\alpha$ is monomial then
\begin{equation}
k_\alpha=R\label{eq:k=R}
\end{equation}
Take a point $\alpha\in\sper\ A$. For an element $a\in\Gamma_\alpha$, let us denote the $\nu_\alpha$-ideal of value $a$ and
$P_{a+}$ the greatest $\nu_\alpha$-ideal properly contained in $P_a$. In other words, we have
\begin{eqnarray}
P_a&=&\{y\in A\ |\ \nu(y)\ge a\},\label{eq:Pa}\\
P_a&=&\{y\in A\ |\ \nu(y)>a\}.\label{eq:Pa+}
\end{eqnarray}
As an application of the above ideas, we give in \S\ref{monomial} a proof of the Connectedness conjecture (and hence also
of the Pierce--Birkhoff conjecture the abstract formulation) in the special case when for all the $\nu_\alpha$-ideals
containing $<\alpha,\beta>$ we have
\begin{equation}
\dim_R\frac{P_a}{P_{a+}}=1\label{eq:dim1}
\end{equation}
(as an $R$-vector space) and $P_a$ is generated by monomials in $x$. In particular, this proves the Connectedness and the
Pierce--Birkhoff conjectures in the case when one of $\alpha$ and $\beta$ is a monomial point of the real spectrum.

The proof goes as follows. Let the notation be as in the statement of the Connectedness conjecture.
By assumption, $g_i\notin<\alpha,\beta>$ so $\nu_\alpha(g_i)<\nu_\alpha(<\alpha,\beta>)$. The assumption (\ref{eq:dim1}) and
the monomiality of $P_a$ for $a<\nu_\alpha(<\alpha,\beta>)$ implies that the $R$-vector space $\frac{P_a}{P_{a+}}$ is
generated by a natural image of a monomial in $x$. Hence, for $1\le i\le s$, $g_i$ can be written as the sum of a dominant
monomial $M(i)$, of value $\nu_\alpha(g_i)<\nu_\alpha(<\alpha,\beta>)$, plus an $R$-linear combination of monomials of value
strictly greater than $\nu_\alpha(g_i)$. Let us denote $M(i)$ by $c_ix^{\gamma(i)}$, with $\gamma(i)\in\N^n$, $c_i\in R$,
and the non-dominant monomials appearing in $g_i$ by $c_{ji}x^{\gamma_{ji}}$, $c_{ji}\in R$, $\gamma_{ji}\in\N^n$, $1\le
j\le N_i$:
$$
g_i=c_ix^{\gamma(i)}+\sum_{j=1}^{N_i}c_{ji}x^{\gamma_{ji}},
$$
where
\begin{equation}
\nu_\alpha\left(x^{\gamma(i)}\right)<\nu_\alpha\left(x^{\gamma_{ji}}\right),\quad1\le j\le N_i.\label{eq:ineqalpha}
\end{equation}
Now, Proposition \ref{sepideal} (2) and the fact that $\nu_\alpha\left(x^{\gamma(i)}\right)<\nu_\alpha(<\alpha,\beta>)$
imply that
\begin{equation}
\nu_\beta\left(x^{\gamma(i)}\right)<\nu_\beta\left(x^{\gamma_{ji}}\right),\quad1\le j\le N_i.\label{eq:ineqbeta}
\end{equation}
By definition of $<\alpha,\beta>$, if $x_q\notin<\alpha,\beta>$ then the sign of $x_q(\alpha)$ is the same as that
of $x_q(\beta)$. Next, (\ref{val2}), (\ref{eq:ineqalpha}) and (\ref{eq:ineqbeta}) imply that
\begin{equation}
\left|x^{\gamma(i)}\right|\ge_\alpha N\left|x^{\gamma_{ji}}\right|\quad\text{for all }N\in R\quad\text{ and }\quad1\le j\le
N_i,\label{eq:ineqalpha1}
\end{equation}
and similarly for $>_\beta$. Let $C$ denote the set of all the points $\delta\in\sper\ A$ such that
\begin{equation}
\left|c_ix^{\gamma(i)}\right|\ge_\delta N_i\left|c_{ji}x^{\gamma_{ji}}\right|\quad\text{for }\quad1\le j\le
N_i,\label{eq:ineqdelta}
\end{equation}
and $x_q(\delta)$ has the same sign as $x_q(\alpha)$ for all $q$ such that $x_q\notin<\alpha,\beta>$. Then $\alpha,\beta\in
C$ and all of the $g_i$ have constant sign on $C$, $1\le i\le s$. In \S\ref{monomial} we give an elementary proof of the
following fact:
\begin{lem}\label{DEUG1} Fix an index $l\in\{0,\dots,n\}$. Let $C$ denote the subset of $\sper\ A$ defined by specifying
$\sgn\ x_q$ (which can be either
strictly positive on all of $C$ or strictly negative on all of $C$) for $q\in\{1,\dots,l\}$ and by imposing, in addition,
finitely many monomial inequalities of the form
\begin{equation}
\left|d_ix^{\lambda_i}\right|\ge\left|x^{\theta_i}\right|,\quad1\le i\le M\label{eq:monineqC1}
\end{equation}
where $d_i\in R\setminus\{0\}$, $\lambda_i,\theta_i\in\N^n$ and $x_q$ may appear only on the right hand side of the
inequalities (\ref{eq:monineqC}) for $q\in\{l+1,\dots,n\}$. Then $C$ is connected.
\end{lem}
This completes the proof of the Connectedness conjecture in this special case. Another connected set
$C$ which can do the same job is the set of all $\delta$ such that the inequalities (\ref{eq:ineqalpha}) are satisfied with
$\alpha$ replaced by $\delta$ and $x_q(\delta)$ has the same sign as $x_q(\alpha)$ for all $q$ such that
$x_q\notin<\alpha,\beta>$. The connectedness of such a set $C$ is an immediate corollary of the main connectedness theorem
of the present paper --- Theorem \ref{con2}, stated below. This will be explained in more detail at the end of the paper.

In \S\ref{affine} we study the behaviour of certain subsets of the
real spectrum under blowing up.

In \S\ref{Des} we recall and adapt to our context some known results
on resolution of singularities, of a purely combinatorial nature. These
results can be considered as a special case of the desingularization
of toric varieties or Hironaka's game \cite{Spi2}. Since they are easy
to prove, we chose to include complete proofs. The conclusion of this
section is an algorithm for resolving singularities of any binomial by
iterating combinatorial (toric) blowings up along non-singular
centers.

Finally, \S\ref{con} is devoted to the proof of the main theorem. Let
$A=R[x_1,\ldots,x_n]$ be a polynomial ring
and let $\omega_{ij},\theta_{il} \in \Q, i \in \{1,\ldots,n\}, j \in
\{1,\ldots,q\}$, $l\in\{1,\dots,u\}$, where $q$ and $u$ are positive integers.

Let $\nu_\delta$ be the valuation associated to the point
$\delta\in\sper(A)$, defined in \S\ref{Cu}.\medskip

\noi\textbf{Notation:} We will write $\nu_\delta(x)$ for the $n$-tuple
$(\nu_\delta(x_1),\ldots,\nu_\delta(x_n))$.
\medskip

Let
$$
h_j(\nu_\delta(x))=\sum\limits_{i=1}^n\omega_{ij}\nu_\delta(x_i)\text{ for }j
\in \{1,\ldots,q\},
$$
and
$$
z_l(\nu_\delta(x))=\sum\limits_{i=1}^n\theta_{il}\nu_\delta(x_i)\text{ for
}l\in\{1,\dots,u\}.
$$
\begin{thm}\label{con2} The set
\begin{equation*}
S=\{\delta\in\sper(A)\ |\ x_i>_\delta 0,i\in\{1,\ldots,n\},\
h_j(\nu_\delta(x))>0,j\in\{1,\ldots,q\},
\end{equation*}
\begin{equation} \label{eq:S*}
z_l(\nu_\delta(x))=0,l\in\{1,\dots,u\}\}
\end{equation}
is connected in the spectral topology.
\end{thm}

In other words, subsets of $\sper\ A$ defined by finitely many
$\Q$-linear equations and {\it strict} inequalities on
$\nu(x_1),\dots,\nu(x_n)$ are connected.
\medskip

To prove Theorem \ref{con2}, we first reduce it to the following Proposition:
\begin{prop}\label{assumeu=01} Assume that $u=0$ in (\ref{eq:S*}), that is, $S$ is defined by finitely many strict $\mathbb
Q$-linear inequalities and no equalities. Then $S$ is connected.
\end{prop}
While revising this paper, we found a very simple proof of Proposition \ref{assumeu=01} using Lemma \ref{DEUG1} (the
first proof of Proposition \ref{assumeu=01} in the current version). However, we decided to keep our original proof because
we feel that it sheds some light on the geometry of the real spectrum. Roughly speaking, this proof (now called the
second proof of Proposition \ref{assumeu=01}) can be thought of as a form of ``path connectedness'' of the set $S$.
Along the way towards the second proof, we exhibit two other types of connected sets in $\sper\ A$, whose connectedness we
believe to be of interest in its own right.

In the forthcoming paper \cite{LMSS}, we will develop the theory of
approximate roots of a valuation. Given a ring $A$ and a valuation
$\nu$, non-negative on $A$, a family of approximate roots is a
collection $\{Q_i\}$, finite or countable, of elements of $A$, having some additional properties as outlined below. A {\bf
generalized monomial} (with respect to a given collection $\{Q_i\}$
of approximate roots) is, by definition, an element of $A$ of the form
$\prod\limits_j Q_j^{\gamma_j}$, $\gamma_j\in\N$. The main defining properties of the
approximate roots are the fact that every $\nu$-ideal $I$ in $A$ is
generated by the generalized monomials contained in it, that is,
generalized monomials $\prod\limits_j Q_j^{\gamma_j}$ satisfying
$$
\sum\limits_j\gamma_j\nu(Q_j)\ge\nu(I),
$$
and the fact that for each $i$, $Q_i$ is described by an explicit formula in terms of
$Q_1,...,Q_{i-1}$. In particular, the valuation $\nu$ is completely
determined by the set $\{Q_i\}$ and the values $\nu(Q_i)$.

We then show that every element $g\in A$ can be written as a finite
sum of the form
\begin{equation}
g=c\mathbf{Q}^\theta+\sum_{j=1}^Nc_j\mathbf{Q}^{\delta_j},
\label{eq:optimalf}
\end{equation}
where $c$ and $c_j$ are elements of $A$ such that $\nu(c)=\nu(c_j)=0$ and $\mathbf{Q}^\theta$ and
$\mathbf{Q}^{\delta_j}$ are generalized monomials such that
\begin{equation}
\nu\left(\mathbf{Q}^\theta\right)<\nu\left(\mathbf{Q}^{\delta_j}\right)\text{ for }1\le j\le N.\label{eq:leading}
\end{equation}

Now let $\delta \in \sper(A)$ and let $\nu=\nu_\delta$. Then, by (\ref{val2}) and (\ref{eq:leading}) the sign of $g$ with
respect to $\leq_\delta$ is determined by the sign of its leading coefficient $c$.
\medskip

Let $\alpha$ and $\beta$ be two points of $\sper(A)$ having a common specialization. In \cite{LMSS} we will give an
explicit description of the set of generalized monomials (with respect to approximate roots $Q_i$ common to $\nu_\alpha$ and
$\nu_\beta$) which generate the separating ideal $<\alpha,\beta>$. Furthermore, we show that all the approximate roots $Q_i$
for $\nu_\alpha$ such that $Q_i \notin <\alpha,\beta>$ are also approximate roots for $\nu_\beta$ and vice versa.
\medskip

Now let $g$ be as in the separation conjecture. The fact that $g\notin<\alpha,\beta>$ implies the existence of an expression
(\ref{eq:optimalf}) in which all the approximate roots $Q_i$ are common for $\nu_\alpha$ and $\nu_\beta$ and the
inequalities (\ref{eq:leading}) hold for both $\nu=\nu_\alpha$ and $\nu=\nu_\beta$.
The inequalities (\ref{eq:leading}) can be viewed as linear inequalities on $\nu(Q_1),\dots,\nu(Q_t)$ with
integer coefficients.

To prove the separation conjecture, we look for connected sets $C \subset \sper(A)$ having the
following properties:

(1) $\alpha, \beta \in C$;

(2) the $Q_i$ appearing in (\ref{eq:optimalf}) are approximate roots simultaneously for all the $\nu_\delta$, $\delta \in C$;

(3) the inequalities (\ref{eq:leading}) hold for $\nu=\nu_\delta$, for all $\delta \in C$;

(4) the leading coefficient $c$ has constant sign on $C$.

Once such a $C$ is found, (\ref{val2}) and (\ref{eq:leading}) imply that the sign of $g$ on $C$ is constant, which proves
that $\alpha$ and $\beta$ lie in the same connected component of $\sper(A) \setminus \{g=0\}$.
\medskip

In \cite{LMSS} we will define the set
$$
C(g,\alpha,\beta)=\left\{\delta\in\sper(A)\ \left|\ Q_i>0,i\in \{1,\ldots,t\},\
\sum_{i=1}^t\omega_{ij}\nu_\delta(Q_i)>0,j\in\{1,\ldots,q\}\right.\right.,
$$
\begin{equation}
\left.\sum_{i=1}^t\lambda_{ij}\nu_\delta(Q_i)=0,j\in\{1,\ldots,l\}\right\},\label{eq:Calphabeta}
\end{equation}
where $\omega_{ij},\lambda_{ij}\in\mathbb Z$, all the
inequalities (\ref{eq:leading}) described above appear among the
$\sum\limits_{i=1}^t\omega_{ij}\nu_\delta(Q_i)>0$ and the remaining equalities and inequalities on the right hand side of
(\ref{eq:Calphabeta}) encode the fact that $Q_1,\dots,Q_t$ are approximate roots for $\nu_\delta$ for all $\delta$ belonging
to $C(g,\alpha,\beta)$. If $g_1,\dots,g_s$ is a finite collection of
elements of $A \ \setminus <\alpha,\beta>$, we put
$$
C(g_1,\dots,g_s,\alpha,\beta)=\bigcap\limits_{i=1}^sC(g_i,\alpha,\beta).
$$
By construction, both sets $C(g,\alpha,\beta)$ and $C(g_1,\dots,g_s,\alpha,\beta)$ contain $\alpha$ and $\beta$, the
element $g$ does not\,\ change\,\ sign on $C(g,\alpha,\beta)$ and none of the elements $g_1,\dots,g_s$ change sign
on\linebreak $C(g_1,\dots,g_s,\alpha,\beta)$. Thus to prove the separation conjecture it is sufficient to prove
that\linebreak $C(g,\alpha,\beta)$ is connected and to prove the Pierce--Birkhoff conjecture it is sufficient to prove that
$C(g_1,\dots,g_s,\alpha,\beta)$ is connected.

A proof of the connectedness of $C(g,\alpha,\beta)$ and $C(g_1,\dots,g_s,\alpha,\beta)$ is relegated to
future papers. In any case, the most delicate part of the proof of the
separation (resp. the Pierce--Birkhoff) conjecture is proving the connectedness of $C(g,\alpha,\beta)$
(resp. $C(g_1,\dots,g_s,\alpha,\beta)$). The connectedness theorem of
the present paper is the special case of the desired result in which
the finite set $\{Q_1,\dots,Q_t\}$ is a subset of the set of variables
$\{x_1,\dots,x_n\}$. In \cite{LMSS} we plan to reduce the connectedness
of $C(g,\alpha,\beta)$ to this special case using sequences of
blowings up of the form described in \S\ref{Des}: we will construct a
sequence $\pi$ of blowings up such that the total transform of each
$Q_i$, $i\in\{1,\dots,t\}$, is a usual monomial with respect to the new
coordinates, times a unit. The preimage of $C(g,\alpha,\beta)$ under
$\pi$ has the form (\ref{eq:S*}) in the new coordinates. This will
reduce the connectedness of $C(g,\alpha,\beta)$ to that of sets of the
form $S$, proved in this paper.
\medskip

We thank the CNRS and the University of Angers for supporting J.
Madden's stay in Angers during a crucial stage of our work on this
paper.
\medskip

We also thank the referee for a very thorough reading of the paper and for
his constructive criticism which helped us improve the exposition and
eliminate numerous inaccuracies.

\section{The Connectedness conjecture implies the Pierce-Birkhoff\\
  conjecture}
\label{implies}

\begin{thm}\label{implies1} The Connectedness conjecture implies the
  Pierce-Birkhoff conjecture.
\end{thm}

\begin{figure}[htbp]
\begin{center}
\input{connfig.pstex_t}
\end{center}
\end{figure}

\noi{\it Proof:} We will assume the Connectedness conjecture and
deduce the Pierce--Birkhoff conjecture in the form of Conjecture
\ref{PBS}. Let $f\in PW(A)$ and let $\{f_i\}_{1\le i\le r}$ denote the
elements of $A$ which represent $f$ on the various closed
constructible subsets $P_i\subset\sper\ A$. Let
$\alpha,\beta\in\sper\ A$ and let
\begin{equation}
T=\{\{i,j\}\subset\{1,\dots,r\}\ |\ f_i-f_j\notin<\alpha,\beta>\}.\label{eq:notsep}
\end{equation}
We apply the Connectedness conjecture to the finite collection
$\left\{f_i-f_j\ \left|\ \{i,j\}\in T\right.\right\}$ of elements of $A$.  By the Connectedness
conjecture, there exists a connected subset $C\subset\sper\ A$ such
that $\alpha,\beta\in C$ and
\begin{equation}
C\cap\{f_i-f_j=0\}=\emptyset\mbox{ for all }\{i,j\}\in
T.\label{eq:conn}
\end{equation}

Let $J\subset\{1,\dots,r\}$ be the set of all the indices $j$ having the following property:
there exist $i_1,\dots,i_s\in\{1,\dots,r\}$ such that
\begin{eqnarray}
\alpha&=&i_1,\label{eq:alphain}\\
j&=&i_s\label{eq:betaend}
\end{eqnarray}
and for each $q\in\{1,...,s-1\}$, we have
\begin{equation}
  C\cap\{f_{i_q}-f_{i_{q+1}}=0\}\ne\emptyset.\label{eq:inter}
\end{equation}
Let $F=\bigcup\limits_{j\in J}(P_j\cap C)$. We have $\alpha\in F$ by definition.

\noi\textit{Claim:} $F=C$; in particular, $\beta\in F$.

\noi\textit{Proof of Claim.} Let $J^c=\{1,\dots,r\}\setminus J$ and $G=\bigcup\limits_{j\in J^c}(P_j\cap C)$. Clearly
\begin{equation}
C=F\cup G\label{eq:FcupG}
\end{equation}
and both sets $F$ and $G$ are closed in the induced topology of $C$ (since all the $P_j$ are closed). Next, we show that
\begin{equation}
F\cap G=\emptyset\label{eq:intvide}
\end{equation}
Indeed, if $\delta\in F\cap G$ then there exist $j\in J$ and $j'\in J^c$ such that $\delta\in P_j\cap P_{j'}$. But
then $f_j(\delta)=f_{j'}(\delta)$, so $\delta\in C\cap\{f_j-f_{j'}=0\}$, hence $j'\in J$, a contradiction. This proves
(\ref{eq:intvide}), so the union in (\ref{eq:FcupG}) is a disjoint union.

Now, since $C$ is connected and $F\ne\emptyset$ (since $\alpha\in F$), the expression (\ref{eq:FcupG}) of $C$ as a
disjoint union of closed sets implies that $G=\emptyset$. Hence $\beta\in F$, which completes the proof of the Claim.
\medskip

Let $j\in J$ be such that $\beta\in P_j$, so that $f_j=f_\beta$. Let $i_1,\dots,i_s$ be as in
(\ref{eq:alphain})--(\ref{eq:inter}), expressing the fact that $j\in J$. Together, (\ref{eq:notsep}), (\ref{eq:conn}) and
(\ref{eq:inter}) imply that $f_{i_q}-f_{i_{q+1}}\in<\alpha,\beta>$ for all
$q\in\{1,...,s-1\}$. In view of (\ref{eq:alphain}) and
(\ref{eq:betaend}), we obtain $f_\alpha-f_\beta\in<\alpha,\beta>$, as
desired.\hfill$\Box$\medskip

\section{The valuation associated to a point in the real spectrum}
\label{Cu}

Let $B$ be a ring and $\alpha$ a point in $\sper\ B$. In this section
we define the valuation $\nu_\alpha$ of $B(\alpha)$, associated to
$\alpha$. We also give a geometric interpretation of points in $\sper\
B$ as semi-curvettes.
\medskip

\noi{\bf Terminology}: If $B$ is an integral domain, the phrase ``valuation
of $B$'' will mean ``a valuation of the field of fractions of $B$, non-negative on $B$''. Also, we will sometimes commit
the following abuse of notation. Given a ring $B$, a prime ideal $\mathfrak p\subset B$, a
valuation $\nu$ of $\frac B{\mathfrak p}$
and an element $x\in B$, we will write $\nu(x)$ instead of $\nu(x\mod\mathfrak p)$, with the usual convention that
$\nu(0)=\infty$, which is taken to be greater than any element of the value group.
\medskip

First, we define the valuation ring $R_\alpha$ by
$$
R_\alpha=\{x\in B(\alpha)\ |\ \exists z\in B[\alpha],|x|\leq_\alpha z\}.
$$
That $R_\alpha$ is, in fact, a valuation ring, follows because for any
$x\in B(\alpha)$, either $x\in R_\alpha$ or $\frac1x\in R_\alpha$. The
maximal ideal of $R_\alpha$ is $M_\alpha=\left\{x\ \in B(\alpha)
  \left|\ |x|<\frac{1}{|z|}, \ \forall z\in B[\alpha] \setminus \{0\}
  \right.\right\}$; its residue field $k_\alpha$ comes equipped
with a total ordering, induced by $\le_\alpha$. For a ring $B$ let $U(B)$ denote the multiplicative group of units of $B$.
Recall that $\dsp{\Gamma_\alpha \cong\frac{B(\alpha)\setminus \{0\}}{U(R_\alpha)}}$ and that the
valuation $\nu_\alpha$ can be identified with the natural homomorphism \\
\centerline{$\dsp{B(\alpha)\setminus \{0\} \to \frac{B(\alpha)\setminus
  \{0\}}{U(R_\alpha)}}$.}
\medskip

By definition, we have a natural ring homomorphism
\begin{equation}
B\rightarrow R_\alpha\label{eq:hom}
\end{equation}
whose kernel is $\mathfrak{p}_\alpha$.  \medskip

\begin{rek} (\cite{Baer}, \cite{Krull}, \cite{BCR} 10.1.10,
p. 217) Conversely, the point $\alpha$ can be reconstructed
from the ring $R_\alpha$ by specifying a certain number of sign
conditions (finitely many conditions when $B$ is noetherian), as we
now explain. Take a prime ideal $\mathfrak{p}\subset B$ and a
valuation $\nu$ of $\kappa({\mathfrak{p}}):=\frac{B_{\mathfrak{p}}}{{\mathfrak{p}}B_{\mathfrak{p}}}$,
with value group $\Gamma$. Let
$$
r=\dim_{\F_2}(\Gamma/2\Gamma)
$$
(if $B$ is not noetherian, it may happen that $r=\infty$). Let
$x_1,\ldots,x_r$ be elements of $\kappa({\mathfrak{p}})$ such that
$\nu(x_1),\ldots,\nu(x_r)$ induce a basis of the $\F_2$-vector space
$\Gamma/2\Gamma$. Then for every $x\in\kappa({\mathfrak{p}})$, there exists $f\in\kappa({\mathfrak{p}})$, and a unit $u$ of
$R_\nu$ such that $x=ux_1^{\epsilon_1}\cdots x_r^{\epsilon_r}f^2$ with $\epsilon_i\in\{0,1\}$
(to see this, note that for a suitable choice of $f$ and $\epsilon_j$ the
value of the quotient $u$ of $x$ by the product $x_1^{\epsilon_1}\cdots
x_r^{\epsilon_r}f^2$ is 0, hence $u$ is invertible in $R_\nu$). Now,
specifying a point $\alpha\in\sper\ B$ supported at ${\mathfrak{p}}$ amounts
to specifying a valuation $\nu$ of $\frac B{\mathfrak{p}}$, a total ordering of the residue
field $k_\nu$ of $R_\nu$, and the sign data $\sgn\ x_1,\dots,\sgn\ x_r$. For $x\notin\mathfrak{p}$, the sign of $x$ is
given by the product $\sgn(x_1)^{\epsilon_1}\cdots\sgn(x_r)^{\epsilon_r}\sgn(u)$, where $\sgn(u)$
is determined by the ordering of $k_\nu$.
\end{rek}

Points of $\sper\ B$ admit the following geometric interpretation (we refer the reader to \cite{Fuc}, \cite{Kap}, \cite{Pre},
p. 89 and \cite{PC} for the construction and properties of generalized power series rings and fields).

\begin{deft} Let $k$ be a field and $\Gamma$ an ordered abelian group. The
  generalized formal power series ring $k\left[\left[t^\Gamma\right]\right]$
  is the ring formed by elements of the form $\sum\limits_\gamma a_\gamma
  t^\gamma$, $a_\gamma\in k$ such that the set $\left\{\gamma\ \left|\ a_\gamma\ne0\right.\right\}$
  is well ordered.
\end{deft}
The ring $k\left[\left[t^\Gamma\right]\right]$ is equipped with the
natural $t$-adic valuation $v$ with values in $\Gamma$, defined by
$v(f)=\inf\{\gamma\ |\ a_\gamma\neq0\}$ for $f=\sum\limits_\gamma
a_\gamma t^\gamma\in k\left[\left[t^\Gamma\right]\right]$. Specifying
a total ordering on $k$ and $\dim_{\F_2}(\Gamma/2\Gamma)$ sign
conditions defines a total ordering on
$k\left[\left[t^\Gamma\right]\right]$. In this ordering $|t|$ is
smaller than any positive element of $k$. For example, if $t^\gamma>0$ for all
$\gamma\in\Gamma$ then $f>0$ if and only if $a_{v(f)}>0$.

For an ordered field $k$, let $\bar k$ denote the real closure of $k$.
The following result is a variation on a theorem of Kaplansky (\cite{Kap}, \cite{Kap2})
for valued fields equipped with a total ordering.

\begin{thm} \textnormal{\textbf{(\cite{PC}, p. 62, Satz 21)}} Let $K$ be a
  real valued field, with residue field $k$ and value group $\Gamma$.
  There exists an injection $K\hookrightarrow\bar
  k\left(\left(t^\Gamma\right)\right)$ of real valued fields.
\end{thm}
\medskip

Let $\alpha\in\sper\ B$ and let $\Gamma$ be the value group of
$\nu_\alpha$. In view of (\ref{eq:hom}) and the Remark above, specifying a
point $\alpha\in\sper\ B$ is equivalent to specifying a total order of $k_\alpha$, a morphism
$$
B[\alpha]\to\bar k_\alpha\left[\left[t^\Gamma\right]\right]
$$
and $\dim_{\F_2}(\Gamma/2\Gamma)$ sign conditions as above.\medskip

\noi We may pass to usual spectra to obtain morphisms
$$
\mbox{Spec}\ \left(\bar k_\alpha\left[\left[t^\Gamma\right]\right]\right)\to\mbox{Spec}\
B[\alpha]\to\mbox{Spec}\ B.
$$
In particular, if $\Gamma=\Z$, we obtain a
\textbf{formal curve} in Spec $B$ (an analytic curve if the series are
convergent). This motivates the following definition:

\begin{deft} Let $k$ be an ordered field. A $k$\textbf{-curvette} on $\sper(B)$ is a morphism of the form
$$
\alpha:B\to k\left[\left[t^\Gamma\right]\right],
$$
where $\Gamma$ is an ordered group. A $k$\textbf{-semi-curvette} is a $k$-curvette $\alpha$ together with
  a choice of the sign data $\sgn\ x_1$,..., $\sgn\ x_r$, where
  $x_1,...,x_r$ are elements of $B$ whose $t$-adic values induce an
  $\F_2$-basis of $\Gamma/2\Gamma$.
\end{deft}

We have thus explained how to associate to a point $\alpha$ of Sper\ $B$ a
$\bar k_\alpha$-semi-curvette. Conversely, given an
ordered field $k$, a $k$-semi-curvette $\alpha$ determines a prime
ideal $\mathfrak{p}_\alpha$ (the ideal of all the elements of $B$
which vanish identically on $\alpha$) and a total ordering on
$B/\mathfrak{p}_\alpha$ induced by the ordering of the ring
$k\left[\left[t^\Gamma\right]\right]$ of formal power series. These
two operations are inverse to each other. This establishes a
one-to-one correspondence between semi-curvettes and points of $\sper\
B$.\medskip

Below, we will often describe points in the real spectrum by specifying the corresponding semi-curvettes.\medskip

\noi{\bf Example:} Consider the curvette $\R[x,y]\to\R[[t]]$
defined by $x\mapsto t^2,\ y\mapsto t^3$, and the semi-curvette given
by declaring, in addition, that $t$ is positive. This semi-curvette is nothing but the upper
branch of the cusp.\medskip

Throughout this paper, we study sets of points of the real spectrum defined by certain properties of their associated
valuations. As the point $\delta$ varies, so does $\nu_\delta$ and its value group $\Gamma_\delta$. In order to describe
such sets in the real spectrum, we need to embed $\Gamma_\delta$ in some ``universal'' ordered group.\medskip

\noi{\bf Notation and convention:} Let us denote by $\Gamma$ the ordered group $\R^n_{lex}$. This means that elements of
$\Gamma$ are compared as words in a dictionary: we say that $(a_1,\dots,a_n)<(a'_1,\dots,a'_n)$ if and only if there
exists $j\in\{1,\dots,n\}$ such that $a_q=a'_q$ for all $q<j$ and $a_j<a'_j$.

The reason for introducing $\Gamma$ is that by
Abhyankar's inequality we have $\mbox{rank}\ \nu_\delta\le\dim A=n$ for all $\delta\in\sper\ A$, so the value group
$\Gamma_\delta$ can be embedded into $\Gamma$ as an ordered subgroup (of course, this embedding is far from being unique).
Let $\Gamma_+$ be the semigroup of non-negative elements of $\Gamma$. \medskip

Fix a strictly positive integer $\ell$. In order to deal rigourously
with $\ell$-tuples of elements of $\Gamma_\delta$ despite the non-uniqueness of the embedding $\Gamma_\delta\subset\Gamma$,
we introduce the category $\mathcal{O}\mathcal{G}\mathcal{M}(\ell)$, as follows. An
object in $\mathcal{O}\mathcal{G}\mathcal{M}(\ell)$ is an ordered
abelian group $G$ together with $\ell$ fixed
generators $a_1,\dots,a_\ell$ (such an object will be denoted by
$(G,a_1,\dots,a_\ell)$). A morphism from $(G,a_1,\dots,a_\ell)$
to $(G',a'_1,\dots,a'_\ell)$ is a homomorphism $G\rightarrow G'$ of
ordered groups which maps $a_j$ to $a'_j$ for each $j$.
\medskip

Given $(G,a_1,\ldots,a_\ell),\ (G',a'_1,\ldots,a'_\ell) \in
Ob(\mathcal{OGM}(\ell))$, the notation
\begin{equation}
(a_1,\ldots,a_\ell) \underset{\circ}{\sim} (a'_1,\ldots,a'_\ell)
\end{equation} will mean that $(G,a_1,\ldots,a_\ell)$ and
$(G',a'_1,\ldots,a'_\ell)$ are isomorphic in $\mathcal{OGM}(\ell)$.
\medskip

Take an element
$$
a=(a_1,\ldots,a_\ell)\in\Gamma_+^\ell.
$$
Let $G\subset\Gamma$ be the ordered group generated by
$a_1,\ldots,a_\ell$. Then $(G,a_1,\ldots,a_\ell) \in
Ob(\mathcal{OGM}(\ell))$. For each $\delta \in \sper(A)$, let
$\Gamma_\delta$ denote the value group of the associated valuation
$\nu_\delta$ and $\Gamma_\delta^*$ the subgroup of $\Gamma_\delta$
generated by $\nu_\delta(x_1)$, \dots, $\nu_\delta(x_n)$. In this way,
we associate to $\delta$ an object $\left(\Gamma_\delta^*,
  \nu_\delta(x_1), \dots, \nu_\delta(x_n)\right)\in
Ob(\mathcal{OGM}(n))$.
\medskip

\section{A proof of the Connectedness conjecture in the case when one of $\alpha$ and $\beta$ is monomial}
\label{monomial}

Consider two points $\alpha,\beta\in\sper\ A$. For an element $a\in\Gamma_\alpha$, let $P_a$ denote the $\nu_\alpha$-ideal
of value $a$ and $P_{a+}$ the greatest $\nu_\alpha$-ideal which is properly contained in $P_a$ (cf.
(\ref{eq:Pa})--(\ref{eq:Pa+})). The purpose of this section is to prove the Connectedness (and hence also the
Pierce--Birkhoff) conjecture assuming that $P_a$ is generated by monomials in $x$ whenever $a\le\nu_\alpha(<\alpha,\beta>)$
and that
\begin{equation}
\dim_R\frac{P_a}{P_{a+}}=1\quad\text{whenever }a<\nu_\alpha(<\alpha,\beta>).\label{eq:dim=1}
\end{equation}
We start with a general (and known) result on semi-algebraic connectedness of constructible sets in the Euclidean space.
We give a proof since we did not find the exact statement we needed in the literature (a slightly weaker form of it is
given in \cite{ABR}, p. 52, Proposition--Definition 62). We state and prove the result in greater generality than needed in
this section, because of the applications we have in mind later in the paper.

Consider the totally ordered set $R_\infty:=\{-\infty\}\cup R\cup\{+\infty\}$, equipped with the order topology.
The map
\begin{eqnarray}
\phi(x)&=&-\frac1{1+x}+1\quad\text{if }x\ge0\\
&=&\frac1{1-x}-1\qquad\text{if }x\le0
\end{eqnarray}
defines a homeomorphism from $R_\infty$ to the interval $[-1,1]$. Let $1 \leq s \leq n-1$. The homeomorphism $\phi$ gives
rise to a homeomorphism
$$
\phi^{n-s}:R_\infty^{n-s}\rightarrow[-1,1]^{n-s}
$$
in the obvious way, that is, by applying the map $\phi$ separately to each component.
A \textbf{closed semi-algebraic} subset of $R^s\times R_\infty^{n-s}$
is, by definition, the preimage of a (relatively) closed semi-algebraic subset of
$R^s\times[-1,1]^{n-s}$ with respect to the product topology of
$R^s\times[-1,1]^{n-s}$, under the map $(id,\phi^{n-s})$. Let $G$ be a semi-algebraic subset of $R^s$.
A mapping
$$
G\rightarrow R_\infty^{n-s}
$$
is said to be \textbf{semi-algebraic} if its graph is a closed semi-algebraic set.

Take a semi-algebraic set $G\subset R^s$. Fix two continuous semi-algebraic mappings
$$
f,g:G\rightarrow R_\infty^{n-s},
$$
having the following property. Write $f=(f_{j+1},\dots,f_n)$, $g=(g_{j+1},\dots,g_n)$, where $f_j$ and
$g_j$ are continuous semi-algebraic $R_\infty$-valued functions on $G$. We require, for
each $(a_1,\dots,a_s)\in R^s$ and each $j\in\{s+1,\dots,n\}$, the inequality
$$
f_j(a_1,\dots,a_s)\le g_j(a_1,\dots,a_s).
$$
Let $D_{f,g}$ denote the semi-algebraic subset of $R^s\times R^{n-s}$ defined by
\begin{equation}
D_{f,g}=\left\{\left.(a_1,\dots,a_n)\in R^n\ \right|\ (a_1,\dots,a_s)\in G,
f_j(a_1,\dots,a_s)\le a_j\le g_j(a_1,\dots,a_s),s<j\le n\right\}.\label{eq:Dfg}
\end{equation}
Let $p:R^s\times R^{n-s}_\infty\rightarrow R^s$ be the canonical projection and
\begin{equation}
p_{f,g}:D_{f,g}\to G\label{eq:pfg}
\end{equation}
the restriction of $p$ to $D_{f,g}$.
\begin{prop}\label{Dfgcon} Assume that $G$ is semi-algebraically connected. Then so is the semi-algebraic set $D_{f,g}$.
\end{prop}
\noi\textbf{Proof:} Composing both $f$ and $g$ with the continuous semi-algebraic mapping $\phi^{n-s}$, we may assume that
the images of $f$ and $g$ lie in $R^{n-s}$ rather than in $R_\infty^{n-s}$.

The map $p_{f,g}$ has semi-algebraically connected fibers (in fact, each fiber of
$p_{f,g}$ is the parallelepiped defined by the inequalities
$$
f_j(a_1,\dots,a_s)\le x_j\le g_j(a_1,\dots,a_s),\quad j\in\{s+1,\dots,n\},
$$
which is semi-algebraically connected (\cite{BCR}, Proposition
2.4.3, p. 31)). Suppose $D_{f,g}$ were not semi-algebraically connected. Then there would exist open semi-algebraic sets
$U,V\subset R^n$ such that
$$
D_{f,g} = (U\cap D_{f,g}) \coprod (V\cap D_{f,g}),
$$
with $U\cap D_{f,g}\ne\emptyset$, $V\cap D_{f,g}\ne\emptyset$. Let $U_0=p_{f,g}(U\cap D_{f,g})$,
$V_0=p_{f,g}(V\cap D_{f,g})$.
\medskip

\noi{\bf Claim:} We have $U_0 \cap V_0= \void$.

\noi{\it Proof of the Claim:} We prove the Claim by
contradiction. Take a point $\xi \in U_0 \cap V_0$. Then
$p_{f,g}^{-1}(\xi)= \left(p_{f,g}^{-1}(\xi) \cap U \right) \coprod
\left(p_{f,g}^{-1}(\xi)\cap V\right)$. We have expressed the fiber
$p_{f,g}^{-1}(\xi)$ as a disjoint union of two non-empty, relatively open semi-algebraic sets
which contradicts the semi-algebraic connectedness of $p_{f,g}^{-1}(\xi)$.
The Claim is proved.
\medskip

It follows from the Claim that
\begin{eqnarray}
U\cap D_{f,g}&=&p_{f,g}^{-1}(U_0)\quad\text{and}\label{eq:UU0}\\
V\cap D_{f,g}&=&p_{f,g}^{-1}(V_0).\label{eq:VV0}
\end{eqnarray}
We have
\begin{equation}
G = U_0 \coprod V_0.\label{eq:notconnected}
\end{equation}
Moreover, (\ref{eq:UU0})--(\ref{eq:VV0}) imply that $U_0$ is the preimage of the open
semi-algebraic set $U$ under the continuous semi-algebraic mapping
$$
\left(id,\frac{f+g}2\right):G\rightarrow R^n=R^s\times R^{n-s},
$$
hence it is relatively open in $G$, and similarly for $V_0$. Thus (\ref{eq:notconnected}) expresses $G$ as a disjoint
union of non-empty relatively open semi-algebraic sets. Then $G$ is not semi-algebraically connected, which gives the
desired contradiction.
\hfill$\Box$
\begin{rek}\label{strict} Proposition \ref{Dfgcon} holds just as well if we replace one or both non-strict inequalities in
(\ref{eq:Dfg}) by strict ones (the same proof applies verbatim).
\end{rek}
\begin{thm}\label{conmonomial} Let the notation and assumptions be as in the beginning of this section.
Take a finite collection $g_1,\dots,g_s$ of elements of $A\setminus<\alpha,\beta>$. Then there exists a connected set
  $C\subset\sper\ A$ such that $\alpha,\beta\in C$ and
  $C\cap\{g_i=0\}=\emptyset$ for $i\in\{1,\dots,s\}$ (in other words,
  $\alpha$ and $\beta$ belong to the same connected component of the
  set $\sper\ A\setminus \{g_1\dots g_s=0\}$).
\end{thm}
\begin{rek} Of course, the conclusion of Theorem \ref{conmonomial} remains valid if its hypotheses are satisfied after
applying an $R$-automorphism of $A$.
\end{rek}
\begin{cor}\label{pbmonomial} Assume that that at least one of $\alpha$, $\beta$ is monomial. Then the conclusion of the
Connectedness conjecture (and hence also that of the Pierce--Birkhoff conjecture in the abstract formulation) holds for the
pair $\alpha$, $\beta$.
\end{cor}
\noi\textbf{Proof of Corollary \ref{pbmonomial} assuming Theorem \ref{conmonomial}:} In general, for each
$a\in\Gamma_\alpha$ we have a non-canonical injection
\begin{equation}
\frac{P_a}{P_{a+}}\subset k_\alpha\label{grdanskalpha}
\end{equation}
(of $\frac A{P_{0+}}$-vector spaces), defined by picking an element $y\in P_a\setminus P_{a+}$ and for each $z\in P_a$
sending $z\mod\ P_{a+}$ to the natural image of $\frac zy$ in $k_\alpha$. Without loss of generality, we may assume that
$\alpha$ is monomial. Then (\ref{eq:k=R}), (\ref{grdanskalpha}) and
the fact that $R$ is a subfield of $\frac A{P_{0+}}$ imply (\ref{eq:dim=1}) for all $a\in\Gamma_\alpha$, in particular,
for $a<\nu_\alpha(<\alpha,\beta>)$. Also, it is easy to see from the definition of monomial points of the real spectrum
that the ideal $P_a$ is monomial for all $a\in\Gamma_\alpha$, in particular for $a\le\nu_\alpha(<\alpha,\beta>)$. Thus
the hypotheses of Theorem \ref{conmonomial} are satisfied whenever $\alpha$ is monomial. This proves the Corollary.
\hfill$\Box$\medskip

\noi\textbf{Proof of Theorem \ref{conmonomial}:} We have $g_i\notin<\alpha,\beta>$ by assumption, so
$\nu_\alpha(g_i)<\nu_\alpha(<\alpha,\beta>)$ (since $<\alpha,\beta>$ is a $\nu_\alpha$-ideal by Proposition
\ref{sepideal} (1)). The assumption (\ref{eq:dim=1}) and the monomiality of $P_a$ for $a\le\nu_\alpha(<\alpha,\beta>)$ imply
that the $R$-vector space $\frac{P_a}{P_{a+}}$ is generated by a natural image of a monomial in $x$. Hence, for $1\le i\le
s$, $g_i$ can be written as the sum of a dominant monomial (which we will denote by $c_ix^{\gamma(i)}$, $\gamma(i)\in\N^n$,
$c_i\in R$), of value $\nu_\alpha(g_i)<\nu_\alpha(<\alpha,\beta>)$, plus an $R$-linear combination of monomials of value
strictly greater than $\nu_\alpha(g_i)$. Let us denote these non-dominant monomials appearing in $g_i$ by
$c_{ji}x^{\gamma_{ji}}$, $c_{ji}\in R$, $\gamma_{ji}\in\N^n$, $1\le j\le N_i$:
$$
g_i=c_ix^{\gamma(i)}+\sum_{j=1}^{N_i}c_{ji}x^{\gamma_{ji}},
$$
where
\begin{equation}
\nu_\alpha\left(x^{\gamma(i)}\right)<\nu_\alpha(x^{\gamma_{ji}}),\quad1\le j\le N_i.\label{eq:ineqalpha2}
\end{equation}
Proposition \ref{sepideal} (2) and the fact that $\nu_\alpha\left(x^{\gamma(i)}\right)<\nu_\alpha(<\alpha,\beta>)$
imply that
\begin{equation}
\nu_\beta\left(x^{\gamma(i)}\right)<\nu_\beta(x^{\gamma_{ji}}),\quad1\le j\le N_i.\label{eq:ineqbeta1}
\end{equation}
Now, (\ref{val2}), (\ref{eq:ineqalpha2}) and (\ref{eq:ineqbeta1}) imply that
\begin{equation}
\left|x^{\gamma(i)}\right|\ge_\alpha N\left|x^{\gamma_{ji}}\right|\quad\text{for all }N\in R\quad\text{ and }\quad1\le j\le
N_i,\label{eq:ineqalpha3}
\end{equation}
and similarly for $\ge_\beta$. By definition of the separating ideal, if $x_q\notin<\alpha,\beta>$ then the sign of
$x_q(\alpha)$ is the same as that of $x_q(\beta)$. Renumbering the $x_q$, if necessary, we may assume that there exists
$l\in\{0,\dots,n\}$ such that $x_1,\dots,x_l\notin<\alpha,\beta>$ and $x_{l+1},\dots,x_n\in<\alpha,\beta>$. Then for
$q\in\{l+1,\dots,n\}$ the variable $x_q$ appears only on the right side of the inequalities (\ref{eq:ineqalpha2}),
(\ref{eq:ineqbeta1}) and (\ref{eq:ineqalpha3}) (otherwise we would have $\nu_\alpha(x_q)<\nu_\alpha(<\alpha,\beta>)$, a
contradiction).

Let $C$ denote the set of all the points $\delta\in\sper\ A$ such that
\begin{equation}
\left|c_ix^{\gamma(i)}\right|\ge_\delta N_i\left|c_{ji}x^{\gamma_{ji}}\right|\quad\text{for }\quad1\le j\le
N_i,\label{eq:ineqdelta1}
\end{equation}
$x_q(\delta)$ has the same sign as $x_q(\alpha)$ for all $q$ such that $x_q\notin<\alpha,\beta>$, and for
$q\in\{l+1,\dots,n\}$ the variable $x_q$ appears only on the right side of the inequalities (\ref{eq:ineqdelta1}). Then
$\alpha,\beta\in C$ and all of the $g_i$ have constant sign on $C$, $1\le i\le s$. Thus to complete the proof Theorem
\ref{conmonomial}, it remains to prove the following Lemma:
\begin{lem}\label{DEUG} Let $C$ denote the subset of $\sper\ A$ defined by specifying $\sgn\ x_q$ (which can be either
strictly positive on all of $C$ or strictly negative on all of $C$) for $q\in\{1,\dots,l\}$ and by imposing, in addition,
finitely many monomial inequalities of the form
\begin{equation}
\left|d_ix^{\lambda_i}\right|\ge\left|x^{\theta_i}\right|,\quad1\le i\le M\label{eq:monineqC}
\end{equation}
where $d_i\in R\setminus\{0\}$, $\lambda_i,\theta_i\in\Q_+^n$ and $x_q$ may appear only on the right hand side of the
inequalities (\ref{eq:monineqC}) for $q\in\{l+1,\dots,n\}$. Then $C$ is connected.
\end{lem}
\noi\textbf{Convention on rational exponents.} It is convenient to allow the exponents in (\ref{eq:monineqC}) to be
rational. We ascribe the obvious meaning to such an inequality by raising both sides to a suitable integer power to make all
the exponents integer.
\begin{rek} Of course, Lemma \ref{DEUG} holds just as well if we replace the non-strict inequalities in
(\ref{eq:monineqC}) by strict ones (the proof is the same). For the proof of Theorem \ref{conmonomial}, it makes no
difference whether we use strict or non-strict inequalities: both work equally well. We state Lemma \ref{DEUG} as we do
because it will be crucial in the first proof of our main connectedness theorem --- Theorem \ref{con2} --- where non-strict
inequalities are essential because there we need the set $C$ to be \textit{closed} as well as connected.
\end{rek}
\noi\textbf{Proof of Lemma \ref{DEUG}:} We proceed by induction on $n$. If $n=1$, all of the inequalities
(\ref{eq:monineqC}) amount to either $|x_1|\le c$ or $|x_1|\ge d$ for some $c,d\in R$, with the latter type of inequality
allowed only if $l=1$, in which case the sign of $x_1$ on $C$ is fixed. Thus $C$ is either empty or an interval of one of
the forms $\{c'\le x_1\le d'\}$, $\{0<x_1\le d'\}$ or $\{c'\le x_1<0\}$, $c',d'\in R$, which are well known to be connected.

Next, assume that $n>1$ and that the result is known for $n-1$. Without loss of generality, we may assume that there exists
$M_1\in\{0,\dots,M\}$ such that the inequalities (\ref{eq:monineqC}) involve $x_n$ for $1\le i\le M_1$ and do not involve
$x_n$ for $M_1<i\le M$. First, suppose
\begin{equation}
l=n.\label{eq:l=n}
\end{equation}
Let $\bar x=(x_1,\dots,x_{n-1})$ and for each $i$, $1\le i\le M_1$, rewrite the $i$-th inequality (\ref{eq:monineqC}) in one
of the forms
\begin{eqnarray}
|x_n|&\le&\left|d'_i\bar x^{\bar\lambda_i}\right|\quad\text{or}\label{eq:xnless}\\
|x_n|&\ge&\left|d'_i\bar x^{\bar\lambda_i}\right|\label{eq:xnmore}
\end{eqnarray}
where $d'_i\in R$ and $\bar\lambda_i\in\Q^{n-1}$. Without loss of generality, we may assume that there exists
$T\in\{0,\dots,M_1\}$ such that for $i\le T$ the $i$-th inequality is of the form (\ref{eq:xnless}) while for $T<i\le M_1$
the $i$-th inequality is of the form (\ref{eq:xnmore}). Then the inequalities (\ref{eq:monineqC}) for $1\le i\le M_1$ can be
rewritten as
\begin{equation}
\max\limits_{T<i\le M_1}\left\{\left|d'_i\bar x^{\bar\lambda_i}\right|\right\}\le|x_n|\le
\min\limits_{1\le i\le T}\left\{\left|d'_i\bar x^{\bar\lambda_i}\right|\right\}.\label{eq:xnsandwich}
\end{equation}
Here the first inequality in (\ref{eq:xnsandwich}) is taken to be non-existent if $T=M_1$ and the second inequality is
considered non-existent if $T=0$.

By (\ref{eq:l=n}), $x_n$ is either strictly positive on all of $C$ or strictly negative on all of $C$. Without loss of
generality, we may assume that $x_n>0$ on $C$. Then we may replace $|x_n|$ by $x_n$ in the inequalities
(\ref{eq:xnsandwich}). Consider the set $\tilde C_{n-1}\subset\sper\ R[x_1,\dots,x_{n-1}]$, defined by the inequalities
\begin{equation}
\left|d'_i\bar x^{\bar\lambda_i}\right|\le\left|d'_j\bar x^{\bar\lambda_j}\right|,\quad
i\in\{T+1,\dots,M_1\},j\in\{1,\dots,T\},
\end{equation}
all the inequalities (\ref{eq:monineqC}) not involving $x_n$ (that is, those with $M_1<i\le M$),
and the same sign conditions on $\sgn\ x_1$, \dots, $\sgn\ x_{n-1}$ as appear in the definition of $C$. By the induction
assumption, $\tilde C_{n-1}$ is a connected, constructible set of $\sper\ R[x_1,\dots,x_{n-1}]$.

A subset of $R^n$ defined by the conditions of the form $\bar x\in C_{n-1}$,
\begin{equation}
f(\bar x)\le x_n\le g(\bar x),\label{eq:xnsandwich1}
\end{equation}
where $f$ and $g$ are two continuous,
semi-algebraic functions defined on a connected, open semi-algebraic set $C_{n-1}\subset R^{n-1}$ with $f(\bar x)\le g(\bar
x)$, is semi-algebraically connected by Proposition \ref{Dfgcon} (if $T=M_1$ in (\ref{eq:xnsandwich}), we take $f=0$ and the
first inequality in (\ref{eq:xnsandwich1}) becomes strict). Hence $C$ is connected by Proposition \ref{tilde}. This
proves Lemma \ref{DEUG} in the case when $l=n$.

If $l<n$, the above argument also goes through, with the following changes:

(1) Since $l<n$, the inequalities of the form (\ref{eq:xnmore}) do not appear in the definition of $C$, so that $T=M_1$.

(2) The inequality (\ref{eq:xnsandwich}) becomes
\begin{equation}
|x_n|\le\min\limits_{1\le i\le M_1}\left\{\left|d'_i\bar x^{\bar\lambda_i}\right|\right\}.\label{eq:xnsandwich2}
\end{equation}
Note that if $x_q$ appears on the right hand side of (\ref{eq:xnsandwich2}) with a negative exponent then $q<l$ and so
$x_q$ is either strictly positive or strictly negative on all of $C$.

(3) The sets $\tilde C_{n-1}\subset\sper\ R[x_1,\dots,x_{n-1}]$ and $C_{n-1}\subset R^{n-1}$ are defined by all the
inequalities (\ref{eq:monineqC}) not involving $x_n$ (that is, those with $M_1<i\le M$), and sign conditions on $\sgn\
x_1$, \dots, $\sgn\ x_{n-1}$ which appear in the definition of $C$.

(4) The functions $f$ and $g$ of (\ref{eq:xnsandwich1}) are given by
\begin{eqnarray}
f&=&-\min\limits_{1\le i\le M_1}\left\{\left|d'_i\bar x^{\bar\lambda_i}\right|\right\}\quad\text{ and}\\
g&=&\min\limits_{1\le i\le M_1}\left\{\left|d'_i\bar x^{\bar\lambda_i}\right|\right\},
\end{eqnarray}
where if $M_1=0$, we take $f=-\infty$ and $g=+\infty$. This completes the proof of Lemma \ref{DEUG} and Theorem
\ref{conmonomial}.
\hfill$\Box$\medskip

Another set $C$ which satisfies the conclusion of Theorem \ref{conmonomial} is the set of all $\delta\in\sper\ A$ such that
the inequalities (\ref{eq:ineqalpha2}) are satisfied with $\alpha$ replaced by $\delta$ and $x_q(\delta)$ has the same sign
as $x_q(\alpha)$ for all $q$ such that $x_q\notin<\alpha,\beta>$. The connectedness of such a set $C$ is a
corollary of Theorem \ref{con2}, as we will explain at the end of this paper.

\section{Affine monomial blowings up}
\label{affine}

In this section we define one of our main technical tools: affine
monomial blowing up. We describe a large class of situations in which
the valuation $\nu_\alpha$ is preserved under blowing up.
\medskip

\noi{\bf Notation.} For a subset $J\subset\{1,\dots,n\}$, $x_J$ will
stand for the set $\{x_q\ |\ q\in J\}$.

Let $G$ be an ordered group. For an $n$-tuple $a=(a_1,\dots,a_n)\in
G^n$ of elements of $G$, we define
$$
\rr\ a=\dim_{\Q}\sum\limits_{j=1}^n\Q a_j\subset G\otimes_{\Z}\Q.
$$
\medskip

Consider a set $J\subset\{1,\dots,n\}$. Fix an element $j\in J$. Let
\begin{eqnarray}
x'_q=x_q&\qquad\ \ \text{if }q=j\text{ or }q \notin J\label{eq:q=j}\\
=\dsp\frac{x_q}{x_j}&\qquad\text{ if }q\in J\setminus\{j\}\label{eq:qinJ}.
\end{eqnarray}
Let $x=(x_1,\dots,x_n)$, $x'=(x'_1,\dots,x'_n)$ and $A'=R[x']$. We
have a natural ring homomorphism $\pi:A\rightarrow A'$ and the
corresponding map of real spectra $\pi^*:\sper\ A'\rightarrow\sper\ A$.
\medskip

\begin{rek} Since the variables $x$ are monomials in the $x'$
and vice versa, we have
\begin{equation}
\rr(\nu(x_1),\dots,\nu(x_n))=\rr(\nu(x'_1),\dots,\nu(x'_n)).\label{eq:slz}
\end{equation}
\end{rek}
\begin{deft}\label{affblup} The map $\pi$ is called an {\bf affine monomial blowing
up} of $\sper\ A$ (along the ideal $(x_J)$). The choice of the set $J$, or the ideal $(x_J)$, is referred to as {\bf the
choice of the center of blowing-up} and the choice of $j\in J$ ---
as {\bf the choice of a coordinate chart}. Finally, let
$\mathfrak{p}$ be a prime ideal of $A$, not containing any of
$x_1,\dots,x_n$. Let $\nu$ be a valuation of $\frac
A{\mathfrak{p}}$. We say that $\pi$ is an {\bf affine monomial blowing up
with respect to }$\nu$ if
\begin{equation}
\nu(x_j)=\min\{\nu(x_q)\ |\ q\in J\}.\label{eq:jmin}
\end{equation}
Consider an index $q\in\{1,\dots,n\}$. We say that $\pi$ is \textbf{independent of }$x_q$ if $q\notin J$.
\end{deft}

\begin{rek} (1) Let $\mathfrak{p'}$ denote the strict transform of $\mathfrak{p}$ in
$\mbox{Spec}\ A'$. Condition (\ref{eq:jmin}) is equivalent to
saying that
\begin{equation} \label{etlabete} \forall g' \in
  \frac{A'}{\mathfrak{p}'} \mbox{ we have } \nu_\delta(g') \geq 0.
\end{equation}
(2) In the notation of Definition \ref{affblup}, suppose that $\pi$ is independent of $x_q$. Then $x_q=x'_q$. In this
way, we can define the notion of a sequence of affine blowings up, independent of $x_q$, recursively in the length of
the sequence. Namely, if $\pi:A\rightarrow A'$ is a sequence of affine blowings up of length $t$, independent of $x_q$,
$x_q$ equals the $q$-th coordinate of $A'$ and $\pi':A'\rightarrow A''$ is an affine blowing up independent of $x_q$
then we say that $\tilde\pi=\pi'\circ\pi$ is a sequence of affine blowings up, independent of $x_q$.
\end{rek}

We consider the coordinates $x'$ as part of the data of the affine
monomial blowing up $\pi$. An affine monomial blowing up is
completely determined by the choice of $J$ and $j$ as above.
\medskip

We will use the following notation. For a set $E\subset\Gamma^n_+$, let
$$
S_E=\{\delta \in \sper(A)\ |\  x_i>_\delta0, i \in
\{1,\ldots,n\},\exists a=(a_1,\ldots,a_n)\in E\text{ such that}
$$
$$
(\nu_\delta(x_1),\dots,\nu_\delta(x_n)) \underset{\circ}{\sim}
(a_1,\ldots,a_n)\}.
$$
and in particular, for $a\in\Gamma_+^n$ we will write
$$
S_a=\{\delta \in \sper(A)\ |\ x_i >_\delta 0, i \in
\{1,\ldots,n\},(\nu_\delta(x_1), \dots,
\nu_\delta(x_n))\underset{\circ}{\sim}   (a_1,\ldots,a_n) \}
$$
We will need the following comparison result which says that blowing up induces a homeomorphism on sets of the form
$S_E$.
\medskip

Let $E$ be a subset of $\Gamma_+^n$. Take a subset $J \subset \{1,\ldots,n\}$.
Pick a $j\in J$ and consider the affine monomial blowing up
$$
\pi:A\rightarrow A'=R[x'],
$$
determined by $J$ and $j$. Assume that $\pi$ is an affine monomial
blowing up with respect to all $\nu_\delta$, $\delta\in S_E$ (in
other words, $\nu_\delta(x_j)=\min\{\nu_\delta(x_q)\}_{q\in J}$ for
all $\delta\in S_E$; we have $\nu_\delta(x'_i)\ge 0$ for all
$\delta\in S_E$ and $i \in \{1,\ldots,n\}$).

For $a\in\Gamma_+^n$, let $a'$ be the element of $\Gamma_+^n$ defined by
\begin{eqnarray}
a'_q =& a_q\qquad\quad\ \text{ if }q \notin J\text{ or }q=j\label{eq:astays}\\
a'_q =& a_q-a_j\quad\ \,\text{if }q \in J \setminus \{j\}\label{eq:achanges}.
\end{eqnarray}

Let $E'=\{a'\ |\ a\in E\}$ and let $S_{E'}$ denote the corresponding
subset of $\sper\ A'$. Let
$$
\pi^*:\sper\ A'\rightarrow\sper\ A
$$
be the map of real spectra induced by $\pi$. It is well known and easy
to see that $\pi^*$ is a homeomorphism away from the zero set $V(x_j)$
of the ideal $(x_j)$. Since $S_E$ is disjoint from $V(x_j)$, $\pi^*$
induces a homeomorphism
$\pi^*|_{(\pi^*)^{-1}(S_E)}:(\pi^*)^{-1}(S_E)\overset\sim\rightarrow S_E$. For
$\delta\in\sper\ A\setminus V(x_j)$, let $\delta'$ denote the unique
preimage of $\delta$ in $\sper\ A'$.

\begin{prop}\label{valnstays1} Take $a\in\Gamma_+^n$ and
  $\delta\in\sper\ A\setminus V(x_j)$. Assume that $\pi$ is an affine blowing up with respect to $\nu_\delta$. Then
\begin{equation}
\nu_\delta=\nu_{\delta'}\label{eq:valnstays}
\end{equation}
and
\begin{equation}
\delta'\in S_{a'}\iff\delta\in S_a.\label{eq:atoa'}
\end{equation}
\end{prop}
\noi{\bf Proof:} To prove (\ref{eq:valnstays}), it suffices to show that
$R_\delta= R_{\delta'}$. By definition,
 \begin{equation}
 R_\delta =\left\{f \in \kappa(\mathfrak{p}_\delta) \ \left|\ \exists g \in
 \frac{A}{\mathfrak{p}_\delta}, |f| \leq g \right.\right\}
 \end{equation}
\begin{equation}
 R_{\delta'} =\left\{f \in \kappa(\mathfrak{p}_\delta) \ \left|\ \exists g' \in
 \frac{A'}{\mathfrak{p}_{\delta'}}, |f| \leq g' \right.\right\}.
 \end{equation}
Since $\frac{A}{\mathfrak{p}_\delta} \hookrightarrow
\frac{A'}{\mathfrak{p}_{\delta'}}$, we have $R_\delta \subset
R_{\delta'}$. To prove the opposite inclusion, note that, since $\pi$ is an affine blowing up with respect to $\nu_\delta$,
we have $\frac{A'}{\mathfrak{p}_{\delta'}} \subset R_\delta$ by
(\ref{etlabete}). Hence, $\forall g' \in \frac{A'}{\mathfrak{p}_{\delta'}},\
\exists g \in \frac{A}{\mathfrak{p}_\delta}$ such that $|g'| \leq
g$. This proves that $R_{\delta'} \subset R_\delta$ as desired.

(\ref{eq:atoa'}) follows from the equations
(\ref{eq:q=j})--(\ref{eq:qinJ}) and
(\ref{eq:astays})--(\ref{eq:valnstays}).
\hfill$\Box$
\begin{cor}\label{Echanges} Assume that $\pi$ is an affine blowing up with respect to $\nu_\delta$ for all $\delta\in S_E$.
We have $\pi^*(S_{E'})=S_E$ and the restriction of $\pi^*$ to $S_{E'}$ is a homeomorphism.
\end{cor}
\begin{deft} The point $\delta'$ is called the {\bf transform} of
$\delta$. Similarly, $a'$, $E'$ and $S_{E'}$ are called {\bf transforms} of
$a$, $E$ and $S_E$, respectively.
\end{deft}
For future reference, we will also define the transform of a $\Q$-linear
relation on $\nu_\delta(x_1),\dots,\nu_\delta(x_n)$. Consider a $\Q$-linear
equality of the form
\begin{equation}
\sum\limits_{i=1}^n\theta_ia_i=0.\label{eq:lineq2}
\end{equation}
\begin{deft} The transform of (\ref{eq:lineq2}) under $\pi$ is the equality
$\sum\limits_{i=1}^n\theta'_ia'_i=0$, where
\begin{eqnarray}
\theta'_i=&\theta_i\ \quad\ \text{ if }i\ne j\\
=&\sum\limits_{q\in J}\theta_q\ \text{ if }i=j.
\end{eqnarray}
\end{deft}
The transforms of all the above objects under sequences of blowings up are
defined in the obvious way (that is, as iterated transforms) recursively with respect to
the length of the blowing up sequence.
\medskip

\noi{\bf Example:} This example shows that Proposition \ref{valnstays1} (resp.
Corollary \ref{Echanges}) is false without the assumption that $\pi$ is an affine monomial blowing up with respect to
$\nu_\delta$ (resp. with respect to $\nu_\delta$ for all $\delta\in S_E)$. Let
$n=2$ and $a=((1,0),(1,0))\in\Gamma^2$. Consider the point
$\alpha\in\sper\ A$ given by the semi-curvette $\alpha=\left(t^{(1,0)},t^{(0,1)}\right)$. Let $J=\{1,2\}$, $j=1$,
$x'_1=x_1$, $x'_2=\frac{x_2}{x_1}$. Since $(1,0)>(0,1)$ (see the definition of lexicographical in \S\ref{Cu}), the
corresponding blowing up $\pi:A\rightarrow A'$ is not an affine blowing up with respect to
$\nu_\alpha$. We have $\mathfrak{p}_{\alpha'}=(0)$, $R_{\alpha'}=R(x'_2)[x_1]_{(x_1)}$, $x_1>0$, and the total ordering on
the residue field $R(x'_2)$ is given by the inequalities $c<x'_2$ for all $c\in R$ (this information describes $\alpha'$
completely). We have $0<x'_1<c<x'_2$ for any positive constant $c\in R$. Moreover, $\nu_{\alpha'}(x'_2)=(0,0)$,
$\nu_{\alpha'}(x'_1)=(1,0)$, so $\nu_{\alpha'}\ne\nu_\alpha$. We have $\alpha'\in S_{a'}$ but $\alpha\notin
S_a$, so the analogues of (\ref{eq:atoa'}) and Corollary \ref{Echanges} do not
hold without the assumption that $\pi$ is an affine monomial blowing up with respect to
$\nu_\delta$.

\section{Desingularization of binomials by monomial blowings
up}\label{Des}

In this section, we recall and adapt to our context a result from the theory of resolution of singularities---a version of
local uniformization of toric varieties by blowings up along smooth centers, or Hironaka's game.

\begin{prop}\label{hir} Consider two $n$-tuples $\alpha$,
$\gamma\in\mathbb N^n$ and the corresponding monomials $x^\alpha$
and $x^\gamma$. Then there exist finitely many
sequences $\pi_i:A\rightarrow A'_i$, $1\le i\le s$, of affine monomial blowings up, having the following properties:

(1) For each $i\in\{1,\dots,s\}$, one of $x^\alpha$ and $x^\gamma$
divides the other in $A'_i$.

(2) For any prime ideal $\mathfrak{p}$ of $A$ not containing
$x_1,\dots,x_n$, and any valuation $\nu$ of $\frac
A{\mathfrak{p}}$, there exists $i\in\{1,\dots,s\}$ such that
$\pi_i$ is a sequence of affine monomial blowings up with respect to
$\nu$.

(3) Assume that there is an index $q$ such that $\alpha_q=\gamma_q=0$. Then all the blowing up sequences $\pi_i$ above
can be chosen to be independent of $x_q$.
\end{prop}

\begin{rek} (1) Proposition \ref{hir} is a special case of
\cite{Spi2}. We give a proof here since it is much simpler than that
of \cite{Spi2}.

(2) Let $\nu$ be a valuation as in Proposition \ref{hir}, and
$i\in\{1,\dots,s\}$ an index satisfying the conclusion (2) of
the Proposition for $\nu$. Then $\nu$ is non-negative on $A'_i$. If
$\nu(x^\alpha)<\nu\left(x^\gamma\right)$, we know {\it which} of
$x^\alpha$, $x^\gamma$ divides the other in $A'_i$, namely, $x^\alpha\
\left|\ x^\gamma\right.$ (and not the other way around).
\end{rek}

\noi{\bf Proof of Proposition \ref{hir}}: (1) and (2) We will define a numerical
character $\tau(\alpha,\gamma)$, consisting of a pair of non-negative
integers, associated to the {\it unordered} pair $(\alpha,\gamma)$. If one of $x^\alpha$, $x^\gamma$ divides the other,
there is nothing to prove. Assume that neither of $x^\alpha$, $x^\gamma$
divides the other. We will describe a subset $J\subset\{1,\dots,n\}$,
such that for any choice of $j\in J$ and the corresponding affine monomial blowing up $A\rightarrow A'$,
writing $x^\alpha={x'}^{\alpha'}$, $x^\gamma={x'}^{\gamma'}$ in the new coordinates, we have
\begin{equation}
\tau\left(\alpha',\gamma'\right)<\tau(\alpha,\gamma)\label{eq:tau}
\end{equation}
in the lexicographical ordering. Since for any prime ideal $\mathfrak{p}$ of $A$ not containing
$x_1,\dots,x_n$, and any valuation $\nu$ of $\frac A{\mathfrak{p}}$, there exists $j\in J$ such that  the corresponding
affine monomial blowing up $A\rightarrow A'$ is a blowing up with respect to $\nu$ (see (\ref{eq:jmin})), Proposition
\ref{hir} (1) and (2) will follow immediately by iterating this procedure.

We start by defining the numerical character $\tau(\alpha,\gamma)$. Let
\begin{eqnarray}
  \alpha=&(\alpha_1,\dots,\alpha_n)\text{ and}\\
  \gamma=&(\gamma_1,\dots,\gamma_n).\qquad
\end{eqnarray}
Let $\delta_q=\min\{\alpha_q,\gamma_q\}$, $1\le q\le n$; let
$\delta:=(\delta_1,\dots,\delta_n)$. Let
$\tilde\alpha:=\alpha-\delta$, $\tilde\gamma=\gamma-\delta$.
\begin{rek}\label{nul} With a view to proving (3) of the Proposition, note that if $\alpha_q=\gamma_q=0$, we have
\begin{equation}
\delta_q=\tilde\alpha_q=\tilde\gamma_q=0.\label{eq:nul}
\end{equation}
\end{rek}

Interchanging $\alpha$ and $\gamma$, if necessary, we may assume that
$|\tilde\alpha|\le|\tilde\gamma|$ (here and below, $|\ |$ stands for
the sum of the components). Put
$$
\tau(\alpha,\gamma):=\left(|\tilde\alpha|,|\tilde\gamma|\right).
$$
Since we are assuming that $x^\alpha\ \not|\ x^\gamma$ in $A$, we have $\tilde\alpha\ne(0,\dots,0)$ (equivalently,
$|\tilde\alpha|>0$). We will now describe a subset
$J\subset\{1,\dots,n\}$, such that for any choice of $j\in J$ and the corresponding affine monomial
blowing up $\pi:A\rightarrow A'$ along $(x_J)$, the inequality
(\ref{eq:tau}) holds.

Write $\tilde\alpha=(\tilde\alpha_1,\dots,\tilde\alpha_n)$,
$\tilde\gamma=(\tilde\gamma_1,\dots,\tilde\gamma_n)$. Renumbering the
variables, we may assume that there exists $a$, $1\le a<n$, such that
$\tilde\alpha_j=0$ for $a<j\le n$ and $\tilde\gamma_j=0$ for $1\le
j\le a$. In other words,
\begin{eqnarray}
  \tilde\alpha&=(\tilde\alpha_1,\dots,\tilde\alpha_a,
  \underbrace{0,\dots,0}_{n-a\text  { zeroes}})\\
  \tilde\gamma&=(\underbrace{0,\dots,0}_{a\text
    { zeroes}},\tilde\gamma_{a+1},\dots,\tilde\gamma_n).
\end{eqnarray}
We may also assume that
\begin{equation}
  \tilde\alpha_i>0\quad\text{for }1\le i\le a.\label{eq:positive}
\end{equation}
Let $J$ denote a minimal subset of $\{1,\dots,n\}$ (in the sense of
inclusion), having the following properties:
\begin{eqnarray}
  &\{1,\dots,a\}\subset J\quad\text{ and}\label{eq:1a}\\
  &\sum\limits_{q\in J}\tilde\gamma_q\ge|\tilde\alpha|\qquad\label{eq:minimal}
\end{eqnarray}
\vskip0.1in
Pick a $j\in J$. Let $\pi:A\rightarrow A'$ be the affine
monomial blowing up along $(x_J)$, associated to $j$ and $J$.

\begin{rek}\label{independent} If $\alpha_q=\gamma_q=0$, (\ref{eq:positive}) implies that $q>a$. Then by the minimality of
$J$ and (\ref{eq:nul}), $\pi$ is independent of $x_q$.
\end{rek}
We will
now write out the monomials $x^\alpha$ and $x^\gamma$ in the new
coordinates and observe that the numerical character $\tau$ has strictly
decreased. Define the non-negative integers $\tilde\alpha'_q$ and
$\tilde\gamma'_q$, $1\le q\le n$, as follows:
\begin{eqnarray}
  \tilde\alpha'_q&=\tilde\alpha_q\qquad\qquad\qquad\quad\text{if }q\ne
  j\label{eq:alphapos}\\
  &=0\qquad\qquad\ \qquad\quad\,\text{if }q=j\label{eq:alpha0}\\
  \tilde\gamma'_q&=\tilde\gamma_q\qquad\qquad\qquad\quad\text{if }q\ne
  j\label{eq:gammastays}\\
  &=\sum\limits_{q\in J}\tilde\gamma_q-|\tilde\alpha|\qquad\quad\text{if
  }q=j.\label{eq:newgamma}
\end{eqnarray}
Put $\tilde\alpha'=(\tilde\alpha'_1,\dots,\tilde\alpha'_n)$,
$\tilde\gamma'=(\tilde\gamma'_1,\dots,\tilde\gamma'_n)$. Let $\delta'$
denote the $n$-vector obtained from $\delta$ by adding
$|\tilde\alpha|$ to the $j$-th component; that is, $$
\delta'=(\delta_1,\dots,\delta_{j-1},\delta_j+|\tilde\alpha|,
\delta_{j+1},\dots,\delta_n).$$
With these definitions, we have
\begin{eqnarray}
  x^\alpha&=(x')^{\delta'+\tilde\alpha'}\\
  x^\gamma&=(x')^{\delta'+\tilde\gamma'}.
\end{eqnarray}
Put $\alpha'=\delta'+\tilde\alpha'$, $\gamma'=\delta'+\tilde\gamma'$.
\begin{lem}\label{taudrops} We have
  $\tau(\alpha',\gamma')<\tau(\alpha,\gamma)$ in the lexicographical
  ordering.
\end{lem}
\noi{\bf Proof:} There are two possibilities: either
$j\in\{1,\dots,a\}$ or $j\in\{a+1,\dots,n\}$. If $j\in\{1,\dots,a\}$
then (\ref{eq:positive}), (\ref{eq:alphapos}) and (\ref{eq:alpha0})
imply that
\begin{equation}
  \left|\tilde\alpha'\right|=|\tilde\alpha|-\tilde\alpha_j<|\tilde\alpha|.\label{eq:newalpha}
\end{equation}
Suppose that $j\in\{a+1,\dots,n\}$. Then by (\ref{eq:alphapos}),
\begin{equation}
  \left|\tilde\alpha'\right|=|\tilde\alpha|.\label{eq:alpheq}
\end{equation}
We will prove that $|\tilde\gamma'|<|\tilde\gamma|$. Indeed, by the
minimality of $J$,
\begin{equation}
  \sum_{q\in J\setminus\{j\}}\tilde\gamma_q<|\tilde\alpha|\label{eq:sumless}
\end{equation}
(otherwise we could replace $J$ by $J\setminus\{j\}$). Now, by
(\ref{eq:gammastays}), (\ref{eq:newgamma}) and (\ref{eq:sumless}),
$$
\left|\tilde\gamma'\right| = \sum\limits^n_{\begin{array}{c}q=a+1\\ q\ne
      j\end{array}}\tilde\gamma'_q+\tilde\gamma'_j=
    \sum\limits^n_{\begin{array}{c}q=a+1\\ q\ne j\end{array}}
   \tilde\gamma_q+\sum\limits_{q\in
    J}\tilde\gamma_q-|\tilde\alpha|
$$
\begin{equation} = \sum\limits_{q=a+1}^n\tilde\gamma_q+
  \left(\sum\limits_{q\in
      J\setminus\{j\}}\tilde\gamma_q-|\tilde\alpha|\right)
  <\sum\limits_{q=a+1}^n\tilde\gamma_q = |\tilde\gamma|.\label{eq:gammadrops}
\end{equation}
To summarize, (\ref{eq:newalpha}), (\ref{eq:alpheq}) and
(\ref{eq:gammadrops}) say that
\begin{equation}
  (|\tilde\alpha'|,|\tilde\gamma'|)<(|\tilde\alpha|,|\tilde\gamma|)=
  \tau(\alpha,\gamma)\label{eq:taudrops}
\end{equation}
in the lexicographical ordering. If
$|\tilde\alpha'|\le|\tilde\gamma'|$ then according to our definitions
$$
\tau(\alpha',\gamma')=(|\tilde\alpha'|,|\tilde\gamma'|),
$$
and the Lemma follows from (\ref{eq:taudrops}). If
$|\tilde\gamma'|<|\tilde\alpha'|$ then
$\tau(\alpha',\gamma')=(|\tilde\gamma'|,|\tilde\alpha'|)<(|\tilde\alpha'|,
|\tilde\gamma'|)$, and, again, the Lemma follows from
(\ref{eq:taudrops}). This completes the proof of Lemma
\ref{taudrops}.\hfill$\Box$\medskip

For as long as neither of $x^\alpha$, $x^\gamma$ divides the other, we can iterate
the construction of Lemma \ref{taudrops}. Since $\tau$ cannot decrease
indefinitely, this process must stop after finitely many
steps. Therefore after a finite number of steps we will arrive at the
situation when one of ${x'}^{\alpha'}$, ${x'}^{\gamma'}$ divides the
other. In other words, one of $x^\alpha$, $x^\gamma$ divides the other
in $A'$.

Of course, the above construction is not unique: at each step we made
an arbitrary choice of a coordinate chart. The
inequality (\ref{eq:tau}) and hence the final conclusion that one of
$x^\alpha$, $x^\gamma$ divides the other hold for all
the possible choices of $j$. Now let $\{\pi_i:A\rightarrow
A'_i\}_{1\le i\le s}$ be the totality of all the blowing up
sequences constructed above for all the possible choices of coordinate
charts, such that one of $x^\alpha$, $x^\gamma$
divides the other in $A'_i$ for all $i\in\{1,\dots,s\}$ (since the
number of choices of coordinate charts is finite at each
step and each sequence stops after finitely many steps, the overall
set is finite). Since for each valuation $\nu$ there
always exists a choice of coordinate chart satisfying (\ref{eq:jmin})
at each step, the set $\{\pi_i:A\rightarrow A'_i\}_{1\le i\le s}$
satisfies the conclusion of the Proposition. (3) of the Proposition holds by construction according to Remarks \ref{nul}
and \ref{independent}. This completes the proof of Proposition \ref{hir}.
\hfill$\Box$
\begin{rek}\label{Jglobal} Keep the above notation. Take an index $i\in\{1,\dots,s\}$ and consider the blowing up
sequence
\begin{equation}
\pi_i:A\rightarrow A'_i,\label{eq:pii}
\end{equation}
as above. Let $\pi:B=R[y_1,\dots,y_n]\rightarrow B'=R[y'_1,\dots,y'_n]$ be one of the affine blowings up appearing in the
blowing up sequence (\ref{eq:pii}). Let $J$
be the subset of $\{1,\dots,n\}$ such that $\pi$ is the blowing up along the ideal $(y_J)$. The definition of $J$ in the
proof of Proposition \ref{hir} depends only on the monomials $x^\alpha$ and $x^\gamma$, written in the $y$ coordinates, and
makes no mention of any valuation whatsoever. Now, take a prime ideal $\mathfrak{p}$ of $A$ not containing $x_1,\dots,x_n$,
and let $\nu$ be a valuation of $\frac A{\mathfrak{p}}$. Let $\mathfrak{p}_B$ be the strict transform of $\mathfrak{p}$ in
$B$. Since, by construction, $y_1,\dots,y_n$ are monomials in $x_1,\dots,x_n$ (with integer exponents), we have
$y_1,\dots,y_n\notin\mathfrak{p}_B$, so $\nu$ is also a valuation of $\frac B{\mathfrak{p}_B}$. Once the set $J$ is defined,
the choice of the coordinate chart $B'$ in such a way that $\pi$ is an affine blowing up with respect to $\nu$ depends on
$\nu$. Proposition \ref{hir} (2) asserts that such a choice of coordinate chart is always possible at every step of the
definition of the sequence $\pi_i$. Namely, we can always choose an index $j\in J$ satisfying
$\nu(y_j)=\min\limits_{q\in J}\{\nu(y_q)\}$.
\end{rek}
Fix a positive integer $u\le n$ and rational numbers $\theta_{jl}$, $l\in\{1,\dots,u\}$, $j\in\{1,\dots,n\}$, such that the
matrix $(\theta_{jl})$ has rank $u$. Let $r:=n-u$. Let $S$ be a subset of $\sper\ A$, such that for all the $\delta\in S$
we have $x_j>_\delta0$, $j\in\{1,\dots,n\}$, and
\begin{equation}
\sum\limits_{j=1}^n\theta_{jl}\nu_\delta(x_j)=0,\qquad l\in\{1,\dots,u\}.
\label{eq:lineq}
\end{equation}
\begin{cor}\label{hir1} There exist finitely many
sequences $\pi_i:A\rightarrow A'_i=R[x_{i1},\dots,x_{in}]$, $1\le i\le s$, of
affine monomial blowings up having the following property. For each $\delta\in S$ there exists $i\in\{1,\dots,s\}$ such that
the following conditions hold:

(1)
\begin{equation}
\nu_\delta\left(x_{i,r+1}\right)=\dots=\nu_\delta\left(x_{in}\right)=0.
\label{eq:value0}
\end{equation}

(2) $\pi_i$ is a sequence of affine monomial blowings up with respect to $\nu_\delta$.

Suppose there is an index $q\in\{1,\dots,n\}$ such that $\theta_{ql}=0$ for all $l\in\{1,\dots,u\}$. Then all the
affine blowing up sequences $\pi_i$ can be chosen to be independent of $x_q$.
\end{cor}
\noi{\bf Proof:} We proceed by induction on $u$. If $u=0$, we have $r=n$ and there is nothing to prove. Assume that
$u>0$, so $r<n$. Consider the non-trivial $\mathbb Q$-linear relation
\begin{equation}
\sum\limits_{j=1}^n\theta_{j1}\nu_\delta(x_j)=0,\label{eq:lineq1}
\end{equation}
satisfied by all the points $\delta\in S$. For each $j\in\{1,\dots,n\}$ let
\begin{eqnarray}
\alpha_j&=&\max\{\theta_{j1},0\}\\
\gamma_j&=&\max\{-\theta_{j1},0\}.
\end{eqnarray}
Put $\alpha=(\alpha_1,\dots,\alpha_n)$; $\gamma=(\gamma_1,\dots,\gamma_n)$. The equation (\ref{eq:lineq1}) says that
$\nu_\delta\left(x^\alpha\right)=\nu_\delta(x^\gamma)$ for all $\delta\in S$. Apply Proposition \ref{hir} to the pair of
monomials $x^\alpha$, $x^\gamma$. Take a point $\delta\in S$ and let $A\rightarrow A'=R[x'_1,\dots,x'_n]$ be a sequence of
affine monomial blowings up with respect to $\nu_\delta$, appearing in the conclusion of Proposition \ref{hir} for
$x^\alpha$ and $x^\gamma$. We have either $\frac{x^\alpha}{x^\gamma}\in A'$ or $\frac{x^\gamma}{x^\alpha}\in
A'$. The fact that $\pi$ is a sequence of blowings up with respect to $\nu_\delta$ means that
\begin{equation}
\nu_\delta(x'_j)\ge0\qquad\text{for }1\le j\le n.\label{eq:posvalue}
\end{equation}
Since $\frac{x^\alpha}{x^\gamma}$ is a monomial in $x'_1,\dots,x'_n$ and
since
$$
\nu_\delta\left(\frac{x^\alpha}{x^\gamma}\right)=
\nu_\delta\left(\frac{x^\gamma}{x^\alpha}\right)=0,
$$
(\ref{eq:posvalue}) implies that one of the $\nu_\delta\left(x'_q\right)$ is equal to zero. Renumbering
the $x'_q$, if necessary, we may assume that $\nu_\delta(x'_n)=0$. Moreover,
$\nu_\delta\left(x'_1\right),\dots,\nu_\delta\left(x'_{n-1}\right)$ satisfy
$u-1$ linearly independent $\Q$-linear relations (namely, the
transforms of the relations (\ref{eq:lineq})). Now let
$$
\pi_i:A\rightarrow A'_i=R[x_{i1},\dots,x_{in}],\qquad1\le i\le s,
$$
be the totality of affine monomial blowing up sequences satisfying the conclusion of Proposition \ref{hir}, applied to the
pair of monomials $x^\alpha$, $x^\gamma$ above. For $i$, $1\le i\le s$, let $S_i$ denote the set of all the points
$\delta\in S$ such that $\pi_i$ is a sequence of affine monomial blowings up with respect to $\delta$. By Proposition
\ref{hir} (2), we have $S=\bigcup\limits_{i=1}^sS_i$. Let $S'_i$ denote the transform of $S_i$ in $\sper\ A'_i$.

At this point, we have constructed finitely many affine monomial blowing up sequences $\pi_i$, $1\le i\le s$, having the
following properties:

(1) For each $\delta\in S$ there exists $i\in\{1,\dots,s\}$ such that $\delta\in S_i$.

(2) For each $i\in\{1,\dots,s\}$ and $\delta\in S_i$ we have
\begin{equation}
\nu_\delta(x_{in})=0\label{eq:last0}
\end{equation}
and
$\nu_\delta\left(x_{i1}\right),\dots,\nu_\delta\left(x_{i,n-1}\right)$ satisfy
$u-1$ linearly independent $\Q$-linear relations
\begin{equation}
\sum\limits_{j=1}^{n-1}\theta^{(i)}_{jl}\nu_\delta(x_{in})=0\label{eq:newrelations}
\end{equation}
(namely, the transforms of the relations (\ref{eq:lineq}), from which we
may delete $\nu_\delta(x_{in})$ in view of (\ref{eq:last0})).

Next, for each $i\in\{1,\dots,s\}$, apply the induction hypothesis to $A'_i$. Since $\nu_\delta(x_{in})$ does not appear
in the relations (\ref{eq:newrelations}), all the subsequent blowings up are independent of $x_q$. In particular, the
equality (\ref{eq:last0}) is preserved by all the subsequent affine monomial blowings up with respect to $\delta$. This
completes the proof of the Corollary.
\hfill$\Box$
\begin{rek}\label{Jglobal1} Let the notation be as in Corollary \ref{hir1}. Take an index $i\in\{1,\dots,s\}$ and consider
the blowing up sequence
\begin{equation}
\pi_i:A\rightarrow A'_i,\label{eq:pii1}
\end{equation}
as in the Corollary. Let $B=R[y_1,\dots,y_n]$ be a polynomial ring appearing at some stage of the blowing up sequence
(\ref{eq:pii1}). This means that $\pi_i$ can be written as a composition of affine monomial blowing up sequences of the
form $A\overset\lambda\rightarrow B\overset{\lambda'}\rightarrow A'_i$. Let $S(B)$ denote the subset of $S$ consisting of
all the points $\delta\in S$ such that $\lambda$ is a sequence of affine monomial blowings up with respect to $\nu_\delta$.
The transforms of the linear relations (\ref{eq:lineq}) in $B$ are well defined for all $\delta\in S(B)$ and are independent
of the choice of $\delta\in S(B)$. When we are at stage $B$ of the construction of Corollary \ref{hir1}, the set $J$ of
the next blowing up is defined in a way independent of the choice of $\delta\in S(B)$; the definition of $J$ depends only
on the transforms of the linear relations (\ref{eq:lineq}) (Remark \ref{Jglobal} and the proof of Corollary \ref{hir1}).
For each $j\in J$, let $\mu_j:B\rightarrow B_j$ be the blowing up determined by the choice of the coordinate chart
corresponding to $j$. For $\delta\in S(B)$, we have $\delta\in S(B_j)$ if and only if $\mu_j$ is an affine blowing up with
respect to $\nu_\delta$. Since for each $\delta\in S(B)$ we can always choose an index $j\in J$ satisfying
$\nu(y_j)=\min\limits_{q\in J}\{\nu(y_q)\}$, we have $S(B)=\bigcup\limits_{j\in J}S(B_j)$.
\end{rek}
\begin{cor}\label{linindep} Let the notation be as in Corollary
\ref{hir1}. Take $\delta\in S$ such that
\begin{equation}
r=\rr(\nu_\delta(x_1),\dots,\nu_\delta(x_n))\label{eq:rr}
\end{equation}
and an index $i\in\{1,\dots,s\}$ such that $\pi_i$ is a sequence of affine monomial blowings up with respect to $\delta$.
Then $\nu_\delta(x_{i1})$, \dots, $\nu_\delta(x_{ir})$ are $\mathbb
Q$-linearly independent.
\end{cor}
\noi{\bf Proof:} This follows immediately from (\ref{eq:slz}) and
(\ref{eq:rr}).

\section{Connectedness theorems for sets in $\sper\ A$.}
\label{con}

As usual, let us denote by $\Gamma$ the ordered group $\R^n_{lex}$ and by
$\Gamma_+$ the semigroup of non-negative elements of $\Gamma$. Let
$\omega_{ij},\theta_{il}\in\Q,i\in\{1,\ldots,n\},j\in\{1,\ldots,q\}$,
$l\in\{1,\dots,u\}$ and consider the subset $E$ of $\Gamma_+^n$ defined by
\begin{equation}
E=\left\{(a_1,\ldots,a_n) \in \Gamma_+^n\ \left|\
\sum_{i=1}^n \omega_{ij}a_i > 0,j \in \{1,\ldots,
q\},\sum\limits_{i=1}^n\theta_{il}a_i=0,l\in\{1,\dots,u\}
\right.\right\}.\label{eq:E}
\end{equation}
\begin{thm}\label{con1} Let $A= R[x_1,\ldots,x_n]$ be a polynomial
ring. Then the set $S_E$ is connected for the spectral topology of
$\sper(A)$.
\end{thm}

A proof of Theorem \ref{con1} will be given after a few lemmas.
\begin{lem}\label{gentop}
Let $X$ be a topological space, normal and compact (not necessarily
Hausdorff) and let $\mathcal{F}$ be a filter of non-empty closed
connected sets. Then the intersection $C= \bigcap\limits_{F \in
  \mathcal{F}}F$ is non-empty, closed and connected.
\end{lem}

\noi{\bf Proof:} That $C$ is non-empty and closed is well known and
easy to see. To prove connectedness, suppose $C=X_1\coprod X_2$,
$X_1,X_2$ closed and non-empty. By normality, there are two open sets $U_1 \supset
X_1,\ U_2 \supset X_2$, $U_1 \cap U_2 = \emptyset$. For any $F
\in\mathcal{F}$, Let $G(F) = F \setminus(U_1 \cup U_2)$. We have
$G(F)\neq\emptyset$ because $F$ is connected. By the compactness of $X$,
we have $\bigcap\limits_{F\in\mathcal{F}}G(F)\neq\emptyset$. Take an
element $x\in\bigcap\limits_{F\in\mathcal{F}}G(F)$. Then
$$
x\in \bigcap\limits_{F \in\mathcal{F}}G(F) \subset \bigcap\limits_{F
  \in \mathcal{F}} F = X_1 \coprod X_2,
$$
but $x \notin U_1 \cup U_2 \supset X_1 \coprod X_2$, which is a
contradiction.
\hfill$\Box$\medskip

\begin{lem}
  Let $B$ be a finitely generated $R$-algebra and $S$ a closed connected subset of $\sper(B)$. Then $S$ is the
  intersection of a filter of closed connected constructible sets.
\end{lem}

\noi{\bf Proof:} Consider the filter $\mathcal{F}$ of all the closed
connected constructible sets containing $S$. The filter $\mathcal{F}$
is not empty because it contains the connected component of $\sper\ B$
containing $S$; such a connected component is constructible by Proposition
\ref{qualite} and Remark \ref{fingen}.

We want to show that $S = \bigcap\limits_{F \in  \mathcal{F}} F$. Suppose that
$S\subsetneqq\bigcap\limits_{F \in \mathcal{F}}F$
and take a point $\alpha \in \bigcap\limits_{F \in
  \mathcal{F}}F\setminus S$.
There exists a basic open set $U \ni \alpha$ such that $U
\cap S= \emptyset$. Then the connected component $C$ of the
complement of $U$, containing $S$, is a constructible closed connected
set containing $S$, hence a member of $\mathcal{F}$. But $\alpha
\notin C$, contradicting $\alpha \in \bigcap\limits_{F \in
  \mathcal{F}}F$.
\hfill $\Box$ \medskip

In the next lemma, we will use the following notation. Let
$$
G=R^s\setminus \ \left\{\prod_{j=1}^sx_j=0\right\}\cong\mbox{Maxr } R[x_1,\ldots,x_s]_{x_1\cdots x_s}.
$$
Let $f,g:G\rightarrow R_\infty^{n-s}$ be two continuous semi-algebraic functions. Let $D_{f,g}$ and $p_{f,g}$ be as in
(\ref{eq:Dfg})--(\ref{eq:pfg}). Let $\tilde D_{f,g}$ be the constructible subset of $\sper\ A$ corresponding to $D_{f,g}$
under the bijection between semi-algebraic subsets of $R^n$ and constructible subsets of $\sper\ A$, described in the
Introduction. Let $\tilde p_{f,g}:\tilde D_{f,g}\rightarrow\tilde G$ be corresponding map of sets in the real spectrum.

\begin{lem}\textnormal{\textbf{(the Projection lemma)}}\label{proj}

If $C_0$ is a closed connected
subset of $\sper\ R[x_1,\ldots,x_s]_{x_1\ldots x_s}$, then the set
\begin{equation}
C_{f,g}:=\tilde p_{f,g}^{-1}(C_0)\label{eq:cn}
\end{equation}
is connected in $\sper\ R[x_1,\ldots,x_n]_{x_1\ldots x_s}$.
\end{lem}
\noi{\bf Proof of Lemma \ref{proj}:} According to the previous Lemma,
we can write $C_0$ as the intersection of a filter of closed
connected constructible sets $\tilde F \in \mathcal{F}$. Take an
$\tilde F\in\mathcal{F}$. Since $\tilde F$ is connected in the real
spectrum, its preimage $F$ in the maximal spectrum
$R^s\setminus\{\prod_{j=1}^sx_j=0\}$ is semi-algebraically connected
(\cite{BCR}, Proposition 7.5.1, p. 130).

By Proposition \ref{Dfgcon} the closed semi-algebraic subset $p_{f,g}^{-1}(F)\subset
D_{f,g}$ is semi-algebraically connected.

Since $p_{f,g}^{-1}(F)$ is semi-algebraic and
semi-algebraically connected, $\tilde p_{f,g}^{-1}(\tilde F)$ is also connected
for the spectral topology (applying Proposition 7.5.1 of \cite{BCR}
once again). Since
$$
C_{f,g}=\bigcap\limits_{\tilde F\in\mathcal F}\tilde p_{f,g}^{-1}(\tilde F),
$$
$C_{f,g}$ is connected by Lemma \ref{gentop}.

\hfill $\Box$
\begin{rek}\label{large} By Remark \ref{strict} the Projection Lemma remains true if one or both
strict inequalities in the definition of $\tilde D_{f,g}$ and $C_{f,g}$ (\ref{eq:Dfg}) are replaced by non-strict ones,
since the above proof applies verbatim also in that case.
\end{rek}
Let the notation be as in Lemma \ref{proj}.
\begin{lem}\label{proj0} Let $C=\{\delta\in\tilde p^{-1}(C_0)\ |\ x_q>0,\nu_\delta(x_q)=0,s<q\le n\}$. Then $C$
is connected.
\end{lem}
\noi{\bf Proof of Lemma \ref{proj0}:}
For each $\delta\in C$, we have
\begin{equation}
\nu_\delta(x_{s+1})=\dots=\nu_\delta(x_n)=0\label{eq:nu=0}
\end{equation}
by definition of $C$. (\ref{eq:nu=0}) is equivalent
to saying that the images of $\frac1{x_{s+1}},\dots,\frac1{x_n}$ in $A(\delta)$ lie in $R_\delta$, that is, that there
exists a polynomial $z\in A$ such that $\frac1{x_q}\le_\delta z$, $s<q\le n$. Take $j\in\{1,\dots,n\}$ such that
$x_j\ge_\delta x_q$, $q\in\{1,\dots,n\}$. Without loss of generality, we may take $z$ to be of the form $z=Nx_j^b$, where
$N\in R$ and $b\in\mathbb N$. For $N\in R$ with $N\ge2$, $b\in\mathbb N$ and $j\in\{1,\dots,n\}$, let $S_{N,b,j}$
denote the set
$$
S_{N,b,j}=\left\{\delta\in C\ \left|\ \frac1{x_q}\le_\delta Nx_j^b,\quad q\in\{s+1,\dots,n\}\right.\right\}.
$$
Thus for each $\delta\in C$ there exist $N\in R$ with $N\ge2$, $b\in\mathbb N$ and
$j\in\{1,\dots,n\}$ such that
$\delta\in S_{N,b,j}$. In other words, $C=\bigcup\limits_{N,b,j}S_{N,b,j}$. Furthermore, take a point $\delta_0\in C_0$ and
let $(x_1(t),\dots,x_s(t))$ denote the curvette representing $\delta_0$. Then the
curvette $(x_1(t),\dots,x_n(t))$ with $x_q(t)=1$, $s<q\le n$, lies in $S_{N,b,j}$ for all the choices of $N,b,j$ as above.
Thus to prove the connectedness of $C$, it is sufficient to prove the connectedness of $S_{N,b,j}$ for all $N,b,j$.

Fix $N$, $b$ and $j$ as above. We will now use the Projection Lemma to prove that $S_{N,b,j}$ is connected. First, suppose
$j\in\{1,\dots,s\}$. Apply the Projection Lemma with
$$
f(x_1,\dots,x_s)=\left(\frac1{Nx_j^b},\dots,\frac1{Nx_j^b}\right)
$$
and $g=(+\infty,\dots,+\infty)$. Let $D^\bullet_{f,g}$ denote the subset of $D_{f,g}$ defined by replacing the
second non-strict inequality in (\ref{eq:Dfg}) by the strict one. We have
$$
S_{N,b,j}=(\tilde p^\bullet_{f,g})^{-1}(C_0),
$$
where $\tilde p^\bullet_{f,g}=\left.\tilde p\right|_{\tilde D^\bullet_{f,g}}$. Then $S_{N,b,j}$ is connected by the
Projection Lemma and Remark \ref{large}.

Next, assume that $j\in\{s+1,\dots,n\}$. First, let
$$
\tilde p_j:\sper\ R[x_1,\dots,x_s,x_j]_{x_1\dots x_s}\rightarrow\sper\ R[x_1,\dots,x_s]_{x_1\dots x_s}
$$
be the natural projection. Let
$$
C_{N,b,j}=\left\{\delta\in\tilde p_j^{-1}(C_0)\left|\ x_j^{b+1}\ge\frac1N\right.\right\}.
$$
Applying the Projection Lemma and Remark \ref{large} to $\tilde p^\bullet_{f,g}$ with $f=N^{-\frac1{b+1}}$, $g=+\infty$, we
see that the set $C_{N,b,j}$ is connected. Moreover, $C_{N,b,j}$ is closed since it is the intersection of the preimage of
the closed set $C_0$ under the continuous map $\tilde p_j$ and the closed set given by the inequality
$x_j^{b+1}\ge\frac1N$.

Next, consider the natural projection
$$
\tilde p:\sper\ R[x_1,\dots,x_n]_{x_1\dots x_sx_j}\longrightarrow\sper\ R[x_1,\dots,x_s,x_j]_{x_1\dots x_sx_j}
$$
We apply the Projection Lemma once again with $s$ replaced by $s+1$,
$$
f(x_1,\dots,x_s,x_j)=\left(\frac1{Nx_j^b},\dots,\frac1{Nx_j^b}\right)
$$
and $g=(+\infty,\dots,+\infty)$. We have
$$
S_{N,b,j}=(\tilde p^\bullet_{f,g})^{-1}(C_{N,b,j}).
$$
Then $S_{N,b,j}$ is connected by the Projection Lemma and Remark \ref{large}. This completes the proof of Lemma
\ref{proj0}. \hfill $\Box$\medskip

We need one more general lemma about closures of subsets of $\sper\ A$ given by finitely many strict $\mathbb Q$-linear
inequalities on $\nu(x_1)$, \dots, $\nu(x_n)$.
\begin{lem} Let $E$ be a subset of $\Gamma_+^n$ given by finitely many strict $\mathbb Q$-linear inequalities as in
(\ref{eq:E}), but with $u=0$ (that is, only strict inequalities and no equalities appear in the definition of $E$). Then
$E$ is closed in $\sper\ A_{x_1\dots x_n}$.
\end{lem}
\noi\textbf{Proof:} Let the notation be as in (\ref{eq:E}) and let $\omega_j=(\omega_{1j},\dots,\omega_{nj})$,
$j \in \{1,\ldots,q\}$. Let the notation $x^{\omega_j}$ stand for $x_1^{\omega_{1j}}\dots x_n^{\omega_{nj}}$. Then we can
write $E$ as
$$
E=\{\delta\in\sper\ A_{x_1\dots x_n}\ |\ x_i>_\delta0,i\in\{1,\dots,n\},\nu_\delta(x^{\omega_j})>0\}.
$$
For $\delta\in\sper\ A$ such that $x_i>_\delta0,i\in\{1,\dots,n\}$, saying that
$\nu_\delta(x^{\omega_j})>0$ is equivalent to saying that $x^{-\omega_j}$ is bounded below by every element of $A$. This is
also equivalent to saying that $x^{-\omega_j}$ is bounded below by every monomial of the form $Nx_i^b$, $N\in R$,
$b\in\mathbb N$, $i\in\{1,\dots,n\}$. For $N\in R$,
$b\in\mathbb N$, $i\in\{1,\dots,n\}$ and $j\in\{1,\dots,q\}$, let
$$
S_{N,b,i,j}=\left\{\delta\in\sper\ A_{x_1\dots x_n}\ \left|\ x_i>_\delta0,i\in\{1,\dots,n\},x^{-\omega_j}\ge
Nx_i^b\right.\right\}.
$$
Clearly, $S_{N,b,i,j}$ is relatively closed in $\sper\ A_{x_1\dots x_n}$. We have
\begin{equation}
S_E=\bigcap\limits_{N,b,i,j}S_{N,b,i,j},\label{eq:SEint}
\end{equation}
hence $S_E$ is also closed. \hfill $\Box$\medskip

The point of the next two lemmas is to reduce Theorem \ref{con1} to the case when $u=0$, that is, to the case when
no $\mathbb Q$-linear equalities appear in the definition of $E$, only strict inequalities.

Pick a subset $J\subset\{1,\dots,n\}$. We want to study blowings up of $\sper\ A$ along the ideal $(x_J)$.
There are $\#J$ possible choices of coordinate
charts, one for each element $j\in J$. For each $j\in J$, let $\pi_{0j}:A\rightarrow A_j$
denote the affine monomial blowing up, defined by $J$ and $j$. Let
$E_{0j}\subset E$ be defined by
$$
E_{0j}=\{a\in E\ |\ \pi_{0j}\text{ is a blowing up with respect to } a\}
$$
and let $E'_{0j}$ denote the transform of $E_{0j}$ under $\pi_{0j}$. By definitions, the set $E'_{0j}$ is defined in $E$ by
imposing the additional inequalities $a_j\le a_q$, $q\in J$. As a subset of $\Gamma_+^n$,
the set $E'_{0j}$ is defined by the linear equalities and strict inequalities, which are the transforms in
$A_j$ of the equalities and inequalities (\ref{eq:E}) (the
transforms of the non-strict inequalities $a_q\ge a_j$, $q\in J$,
which define $E_{0j}$ inside $E$, have the form $a'_q\ge0$; these
inequalities hold automatically for $(a'_1,\dots,a'_n)\in\Gamma_+^n$ and need not be
taken into account).
\begin{lem}\label{Econvex} Take two indices $j,\tilde j\in J$. Assume that both $E_{0j}$ and
$E_{0\tilde j}$ are not empty, in other words, that there exist $a,\tilde a\in E$ such that
\begin{eqnarray}
a_j&=&\min\limits_{q\in J}a_q\quad\text{ and}\label{eq:aless}\\
\tilde a_{\tilde j}&=&\min\limits_{q\in J}\tilde a_q.\label{eq:tildeless}
\end{eqnarray}
Then there exists a chain of indices $\{j_1,\dots,j_s\}\subset\{1,\dots,n\}$ such that $j=j_1$, $\tilde j=j_s$ and
$E_{0j_q}\cap E_{0j_{q+1}}\ne\emptyset$ for all $q\in\{1,\dots,s-1\}$.
\end{lem}
\noi{\bf Proof:} Let $\Gamma_\circ$ denote the subset of $\Gamma_+$ consisting of all the $n$-tuples of real numbers of the
form $(d,\underbrace{0,\dots,0}_{n-1\text{ zeroes}})$, $d\in\mathbb R$. Let $\Gamma^n_\circ$ denote the subset of
$\Gamma^n_+$ consisting of all the elements of the form $(c_1,\dots,c_n)$, $c_j\in\Gamma_\circ$. Let
$E^\circ=E\cap\Gamma^n_\circ$, $E_{0j}^\circ=E_{0j}\cap\Gamma^n_\circ$, etc.

Given an $n$-tuple $\epsilon=(\epsilon_1,\dots,\epsilon_n)$ of real numbers and an element $d=(d_1,\dots,d_n)$ of $\Gamma_+$
(where $d_q\in\mathbb R$), let $d(\epsilon)$ denote the element of $\Gamma_\circ$ defined by
$(d\cdot\epsilon,\underbrace{0,\dots,0}_{n-1\text{ zeroes}})$, where $d\cdot\epsilon=\sum\limits_{q=1}^nd_j\epsilon_j$.
If $c=(c_1,\dots,c_n)$ is an element of $\Gamma^n_+$, with $c_j\in\Gamma_+$, we will write $c(\epsilon)$ for
$(c_1(\epsilon),c_2(\epsilon),\dots,c_n(\epsilon))\in\Gamma^n_\circ$.

Take two indices $l,\tilde l\in J$ and elements $c\in E_{0l}$, $\tilde c\in E_{0l}\cap E_{0\tilde l}$. Since $E_{0l}$
(resp. $E_{0l}\cap E_{0\tilde l}$) is defined by linear equations and inequalities (not necessarily strict), we can choose
real numbers
$$
\epsilon_1\gg\epsilon_2\gg\dots\gg\epsilon_n>0
$$
such that $c(\epsilon)\in E_{0l}$ and $\tilde c(\epsilon)\in E_{0l}\cap E_{0\tilde l}$. Thus
\begin{equation}
E_{0l}\ne\emptyset\iff E^\circ_{0l}\ne\emptyset.\label{eq:nonempty}
\end{equation}
and
\begin{equation}
E_{0l}\cap E_{0\tilde l}\ne\emptyset\iff E^\circ_{0l}\cap E^\circ_{0\tilde l}\ne\emptyset.\label{eq:nonempty1}
\end{equation}
As a topological space, identify $\Gamma_\circ^n$ with $\mathbb R^n$, with the usual Euclidean topology. Since
$E^\circ$ is defined in $\mathbb R^n_+$ by linear equations and inequalities, it is connected in the Euclidean topology.
Now, each $E_{0l}^\circ$ is relatively closed in $E^\circ$. Hence the conclusion of the Lemma holds with $E_q$ replaced by
$E_q^\circ$. In view of (\ref{eq:nonempty})--(\ref{eq:nonempty1}), this completes the proof of the Lemma.
\hfill$\Box$\medskip

Next, we reduce the Theorem to the case when $u=0$, that is, when $E$ is defined by strict inequalities (and no equalities).
Without loss of generality, we may assume that the equalities
\begin{equation}
\sum\limits_{i=1}^n\theta_{il}a_i=0,\quad l\in\{1,\dots,u\}\label{eq:lineq3}
\end{equation}
appearing in the definition of $E$ are linearly independent.
\begin{lem}\label{u=0} If Theorem \ref{con1} is true for $u=0$ (that is, in the case when there are no equalities in the
definition of $E$) then it is true in general.
\end{lem}
\noi{\bf Proof:} Assume Theorem \ref{con1} for $u=0$. Apply Corollary
\ref{hir1} to the set $S_E$. For each $i\in\{1,\dots,s\}$, let
\begin{equation}
E_i=\{a\in E\ |\ \pi_i\text{ is a sequence of blowings up with
  respect to }a\}
\end{equation}
(cf. Remark \ref{onlya}). By Corollary \ref{hir1} (2),
$E=\bigcup\limits_{i=1}^sE_i$. Let $E'_i$ denote the transform of
$E_i$ under $\pi_i$ (by definition of $E_i$, $E'_i$ is well
defined). By construction, $E'_i$ is contained in the set
$\{a_{r+1}=\dots=a_n=0\}$, and is defined in $\Gamma^n_+\cap\{a_{r+1}=\dots=a_n=0\}$ by finitely many strict
$\Q$-linear inequalities --- the transforms of the inequalities
$\sum\limits_{q=1}^n\omega_{qj}a_q>0$ under $\pi_i$. Let
$\iota:\Gamma^r\rightarrow\Gamma^n$ denote the natural inclusion of
the set $\{a_{r+1}=\dots=a_n=0\}$ in $\Gamma^n$. Let
$E^\dag_i=\iota^{-1}(E'_i)$. Let $\tilde p_i$ denote the natural projection $\tilde p_i:\sper\
R[x_{i1},\dots,x_{in}]_{x_{i1}\dots x_{in}}\rightarrow\sper\
R[x_{i1},\dots,x_{ir}]_{x_{i1}\dots x_{ir}}$. The set $S_{E'_i}$ can be written as
$$
S_{E'_i}=\left\{\left.\delta\in\tilde p_i^{-1}\left(E^\dag_i\right)\ \right|\
x_{ij}>_\delta0,\nu_\delta(x_{ij})=0\text{ for }r<j\le n\right\}\rightarrow\sper\
R[x_{i1},\dots,x_{ir}]_{x_{i1}\dots x_{ir}}.
$$
Now, $S_{E^\dag_i}$ is connected by the $u=0$ case of Theorem \ref{con1}. Then
$S_{E'_i}$ is connected by Lemma \ref{proj0}. Hence $S_{E_i}$
is connected by Corollary \ref{Echanges}.

Let $N_i$ denote the number of affine monomial blowings up composing
$\pi_i$ and let
$$
N=\max\limits_{1\le i\le s}\{N_i\}.
$$
To complete the proof of Lemma \ref{u=0}, we will now show that $S_E$
is connected by induction on $N$. Let $J\subset\{1,\dots,n\}$ be such
that the first blowing up in each of the sequences $\pi_i$ is a
blowing up along $J$ (such a $J$ exists according to Remark \ref{Jglobal1}).

The set $E'_{0j}$ is defined in $\Gamma_+^n$ by the linear
equalities and strict inequalities, which are the transforms in
$A_j$ of the equalities and inequalities (\ref{eq:E}) (the
transforms of the non-strict inequalities $a_q\ge a_j$, $q\in J$,
which define $E_{0j}$ inside $E$, have the form $a'_q\ge0$; these
inequalities hold automatically in $\Gamma_+^n$ and need not be
taken into account). Hence, by the induction assumption,
$S_{E'_{0j}}\cong S_{E_{0j}}$ is connected.

Now, Lemma \ref{Econvex} shows that $S_E=\bigcup\limits_{j\in J}S_{E_{0j}}$ is connected, as desired. This completes the
proof of Lemma \ref{u=0}.
\hfill$\Box$\medskip

To complete the proof of Theorem \ref{con1}, it remains to prove
\begin{prop}\label{assumeu=0} Assume that $u=0$ in (\ref{eq:E}), that is, $E$ is defined by finitely many strict $\mathbb
Q$-linear inequalities. Then $S_E$ is connected.
\end{prop}
\noi\textbf{First proof of Proposition \ref{assumeu=0}:} Write $S_E=\bigcap\limits_{N,b,i,j}S_{N,b,i,j}$, as in
(\ref{eq:SEint}). Let
$$
S_{N,b}=\bigcap\limits_{j=1}^q\left\{\left.\delta\in\bigcap\limits_{i=1}^nS_{N,b,i,j}\ \right|
x^{-\omega_j}\ge1\right\}.
$$
For each $N,b$ the set $S_{N,b}$ is defined by finitely many \textit{monomial} inequalities and the sign conditions $x_t>0$,
$1\le t\le n$. Thus $S_{N,b}$ is connected by Lemma \ref{DEUG}. If $N\le N_1$ and $b\le b_1$ then $S_{N,b}\subset
S_{N_1,b_1}$. Then the collection $\{S_{N,b}\ |\ N\in R,N>0,b\in\mathbb N\}$ is a filter of
connected, closed subsets of $\sper\ A_{x_1\dots x_n}$. Then $S_E$ is connected by Lemma \ref{gentop}. This completes the
first proof of Proposition \ref{assumeu=0} and hence also of Theorem \ref{con1}.\hfill$\Box$\medskip

Before giving the second proof of Theorem \ref{con1}, we first prove two connectedness results.
\begin{prop}\textnormal{\textbf{(fiber connectedness)}}\label{fcl} Fix an element
$a=(a_1,\dots,a_n)\in\Gamma_+^n$. The set $S_a$ is connected in
$\sper(A)$.
\end{prop}
\noi{\bf Proof of Proposition \ref{fcl}:} Let $r =\rr\ a$. Take a subset $J
\subset \{1,\ldots,n\}$ and let
$$
x_J = \{x_q\ |\ q \in J\}.
$$
Let $j\in J$ be such that
\begin{equation}
a_j=\min\{a_q\}_{q\in J}\label{eq:amin}
\end{equation}
and consider the affine monomial blowing up
$$
\pi:A\rightarrow A'=R[x'],
$$
determined by $J$ and $j$:

$x'_q=x_q$, $a'_q = a_q$\qquad\ \ if $q \notin J$ or $q=j$

$x'_q= \dsp\frac{x_q}{x_j}$, $a'_q=a_q-a_j$\ \,if $q \in J \setminus \{j\}$.
\begin{rek}\label{onlya} By (\ref{eq:amin}), $\pi$ is an affine monomial blowing up with
respect to $\delta$ for \textit{all} $\delta\in S_a$. In other words, the property of being an affine monomial blowing up
with respect to $\delta$ depends only on $a$, not on the particular choice of $\delta\in S_a$. Here and below, we will
sometimes say that $\pi$ is an affine monomial blowing up with respect to $a$.
\end{rek}
Let $(G',a'_1,\ldots,a'_n)$ be the ordered group generated by the
$a'_i$ (we have $G'=G$) and let $S_{a'}$ be the corresponding subset
of $\sper(A')$. Then
\begin{equation}
S_a \cong S_{a'}\label{eq:a=a'}
\end{equation}
by Proposition \ref{valnstays1}.

We now iterate the above procedure. By Corollaries \ref{hir1} and
\ref{linindep} (where we take $S=S_a$), after a
succession of such transformations we may assume (after passing to the
new coordinates) that $a_1,\ldots,a_r$ are
$\Q$-linearly independent and
\begin{equation}
a_{r+1}= \cdots = a_n=0.\label{eq:ar+1=0}
\end{equation}
Again, note that by Remark \ref{onlya}, (\ref{eq:a=a'}) and
induction on the number of blowings up required, there is a single sequence of blowings up satisfying the conclusions of
Corollary \ref{hir1} which is a sequence of affine monomial blowings up \textit{simultaneously} for all $\delta\in S_a$. As
$a_1,\ldots,a_r$ are $\Q$-linearly independent, there exists a unique
point
$$
\alpha=\left(t^{a_1},\ldots,t^{a_r}\right)\in\sper\ R[x_1,\dots,x_r]_{x_1\dots x_r}
$$
such that
\begin{equation}
x_i>_\alpha0\label{eq:xpositive}
\end{equation}
and $\nu_\alpha(x_i)=a_i$ for
$i\in\{1,\ldots,r\}$: because of the linear independence of the $a_i$, the
support of $\alpha$ is the zero ideal of $R[x_1,\ldots,x_r]$ and the
ordering $\leq_\alpha$ of $R(x_1,\ldots,x_r)$ is completely described
by the inequalities (\ref{eq:xpositive}) and
\begin{equation} \label{eq:linindep}
cx^\gamma < x^\epsilon\quad\text{for all positive }c\in R \Longleftrightarrow\sum\gamma_ja_j > \sum \epsilon_j a_j.
\end{equation}
The fact that $a_1,\ldots,a_r$ are $\Q$-linearly independent implies
that (\ref{eq:linindep}) imposes a total ordering on the set of
monomials and hence on $R(x_1,\ldots,x_r)$. The point $\alpha$ is
closed in the relative topology of $\sper\ R[x_1,\dots,x_r]_{x_1\dots
  x_r}$ because $\alpha$ has no non-trivial specializations. Now, let $\tilde p$ denote the natural projection
$\tilde p:\sper\ A_{x_1,\dots,x_r}\rightarrow\sper\ R[x_1,\dots,x_r]_{x_1\dots x_r}$. Then
$$
S_a=\left\{\left.\delta\in\tilde p^{-1}(\{\alpha\})\ \right|\ x_q>0,\nu_\delta(x_q)=0\text{ for }r<q\le n\right\},
$$
hence $S_a$ is connected by Lemma \ref{proj0}. This completes the proof of Lemma \ref{fcl}.
\hfill $\Box$ \medskip

\begin{cor}\label{SEdag} Let $S^\dag$ be a connected component of
  $S_E$ and $a$ an element of $\Gamma_+^n$.
Then either $S_a \subset S^\dag$ or $S_a\cap S^\dag=\emptyset$. There
exists a subset $E^\dag \subset E$ such that $S^\dag=S_{E^\dag}$.
\end{cor}

Take $c=(c_1,\dots,c_n),d=(d_1,\dots,d_n) \in E$ such that $c_j=d_j$ for $j\in\{1,\dots,n-1\}$. We define the {\bf segment}
$[c,d]\subset E$ as
\begin{equation}
[c,d]= \{\left(\begin{array}{c} c_1 \\ \vdots \\ c_{n-1} \\ e_n \end{array}
\right)\ |\ c_n\le e_n\le d_n\}.\label{eq:[c,d]}
\end{equation}
Let
\begin{equation}
S_{[c,d]} = \left\{\left.\alpha \in S_{(c_1,\ldots,c_{n-1},e_n)}\
  \right|\ (c_1,\ldots,c_{n-1},e_n) \in [c,d]\right\}.\label{eq:Scd}
\end{equation}
For an element $c=(c_1,\dots,c_n)\in\Gamma^n$, we will write $c_j=(c_{j1},\dots,c_{jn})$ with $c_{ji}\in\mathbb R$. In this
way, we may think of $c$ as an $(n\times n)$ matrix whose $j$-th row is $c_j$.
\begin{prop}\label{cdconnect}
  Take $c,d \in E$ such that :

  (1) $c_j=d_j$ for $j \in \{1,\ldots,n-1\}$

  (2) $c_{n1} < d_{n1}$

  (3) $c_{11},\ldots,c_{n1},d_{n1}$ are $\Q$-linearly independent.

  \noi Then $S_{[c,d]}$ is connected.
\end{prop}

\noi{\bf Proof of Proposition \ref{cdconnect}:} We start with two Lemmas.

Consider an element $b\in\Gamma_+^n$ such that $b_1,\dots,b_n$ are
$\Q$-linearly independent. Then $S_b$ consists of a single point
$\delta(b)$.

\begin{lem}\label{dominant1} Let $U$ be a basic open set for the
  spectral topology of $\sper\ A$, containing $\delta(b)$.
There exists a subset $V\subset\Gamma^n_+$, defined by strict
$\Q$-linear inequalities, such that $b\in V$ and for any $a\in V$, we
have $S_a\subset U$.
\end{lem}

\noi{\bf Proof:} Let $U=\{\alpha\ |\ f_1(\alpha) >0,\ldots,f_s(\alpha)
> 0\}$. Write
\begin{equation} \label{eq:2*} f_j= \sum_{\gamma \in
\N^n}  c_{j\gamma}x^{\gamma}
\end{equation}
Let $M(j)=c_{j,\gamma(j)}x^{\gamma(j)}$ be the monomial of $f_j$ of
smallest valuation (which exists because $b_1,\ldots,b_n$ are
$\Q$-linearly independent).

For any $j,\ 1\leq j \leq s$,
\begin{equation} \label{eq:2**}
  \nu_{\delta(b)}\left(x^{\gamma(j)}\right) <
  \nu_{\delta(b)}\left(x^{\gamma}\right)
\end{equation}
for all $\gamma$ such that $c_{j\gamma} \neq 0$. Writing
$\gamma=(\gamma_1,\dots,\gamma_n)$,
$\gamma(j)=(\gamma_1(j),\dots,\gamma_n(j))$, (\ref{eq:2**}) is
equivalent to saying that $\sum\limits_{q=1}^n (\gamma_q
-\gamma_q(j))b_q > 0$. Let $V$ be the subset of $\Gamma^n_+$, defined
by the inequalities $\sum\limits_{q=1}^n(\gamma_q -\gamma_q(j))a_q >0$
for all $\gamma$ with $c_{j\gamma} \neq 0$. Take an element $a\in V$
and a point $\delta\in S_a$. By construction, each $f_j$ has the same
dominant monomial at $\delta$ as at $\delta(b)$. Then $S_a\subset
U$, as desired.
\hfill$\Box$\medskip

Let $s,r\in\{1,\dots,n-1\}$, $s\le r$. Consider an element $b\in\Gamma_+^n$ such
that
$$
b_{r+1}=\dots=b_n=0
$$
and $b_1,\dots,b_r$ are $\Q$-linearly
independent. For $1\le j\le r$, write $b_j=(b_{j1},\dots,b_{jn})$ with
$b_{jq}\in\R$. Consider the $(r\times n)$-matrix $(b_{jq})^{1\le j\le
r}_{1\le q\le n}$. Assume that the rows of the $(r\times s)$-matrix
  $(b_{jq})^{1\le j\le r}_{1\le q\le s}$ are $\mathbb Q$-linearly independent. Define
  $b^\bullet=(b_1^\bullet,\dots,b_n^\bullet)$ by
$b^\bullet_i=b_i, i \neq r+1$ and $b^\bullet_{r+1}
=(\underbrace{0,\dots,0}_{s},1,0,\ldots,0)$.

\begin{lem}\label{dominant2} We have $\sur{S_b} \cap S_{b^\bullet}
  \ne \emptyset$.
\end{lem}
\noi{\bf Proof:} Pick a point in $S_{b^\bullet}$, for example,
$$
\alpha^\bullet=\left(t^{b_1},\ldots,t^{b_r},t^{b^\bullet_{r+1}},1,\ldots,1\right)
$$
(here we are representing $\alpha^\bullet$ by a parametrized
semi-curvette in $R^n$ - see \S\ref{Cu}). Let
$$
U=\{\alpha\ |\ f_1(\alpha) >0,\ldots , f_q(\alpha) > 0\}
$$
be a basic open set such that $U \ni \alpha^\bullet$. It remains to
prove that $U \cap S_b \neq \emptyset$. Let $\sur{x}=(x_1,\ldots,x_r)$ and $\sur{\sur{x}} =
(x_{r+2},\ldots,x_n)$. Write each $f_j$ as a polynomial in the variables $\sur x$ and $x_{r+1}$, whose coefficients
are polynomials in $\sur{\sur{x}}$:
\begin{equation} \label{*}
f_j= \sum_{\gamma=(\sur{\gamma},\gamma_{r+1})}
c_{j\gamma}(\sur{\sur{x}})\sur{x}^{\sur{\gamma}}x_{r+1}^{\gamma_{r+1}}
\end{equation}
(here $c_{j\gamma}(\sur{\sur{x}})$ is the coefficient in $f_j$ of the monomial
$\sur{x}^{\sur{\gamma}}x_{r+1}^{\gamma_{r+1}}$). Let
$M(j)$ be the term of smallest valuation among those appearing on the right hand side of (\ref{*}); such a term is unique
because $b^\bullet_1,\ldots,b^\bullet_{r+1}$ are $\Q$-linearly independent. Let $\gamma_{r+1}(j)$ (resp. $\sur{\gamma}(j)$)
denote the exponent of $x_{r+1}$ (resp. of $\sur{x}$) in $M(j)$ and let $c_{j,\gamma(j)}\in R[\sur{\sur{x}}]$ denote the
coefficient of $\sur{x}^{\sur{\gamma}(j)}x_{r+1}^{\gamma_{r+1}(j)}$ in (\ref{*}), so that
$M(j)=c_{j,\gamma(j)}\sur{x}^{\sur{\gamma}(j)}x_{r+1}^{\gamma_{r+1}(j)}$.

To say that $M(j)$ is the monomial of $f_j$ of smallest valuation
means that for each $j,\ 1\leq j \leq q$,
\begin{equation} \label{**} \nu\left(\sur{x}^{\sur{\gamma}(j)}\right) \leq
\nu\left(\sur{x}^{\sur{\gamma}}\right)
\end{equation}
for all $\gamma=(\sur{\gamma},\gamma_{r+1})$ such that $c_{j\gamma}
\neq 0$ and if equality holds in (\ref{**}) then
\begin{equation} \label{***}
\gamma_{r+1}(j) < \gamma_{r+1}.
\end{equation}
Take $\epsilon\in R$ sufficiently small so that for all
$j\in\{1,\dots,q\}$ we have
$$
\left|c_{j,\gamma(j)}(1,\ldots,1) \epsilon^{\gamma_{r+1}(j)}\right| >
\sum_{\begin{array}{c}\gamma=(\sur{\gamma},\gamma_{r+1})\\ \gamma_{r+1} >
\gamma_{r+1}(j)\end{array}} \left|c_{j,(\sur{\gamma}(j),\gamma_{r+1})}(1,\ldots,1)
\epsilon^{\gamma_{r+1}}\right|.
$$
Now let $\alpha=(t^{b_1},\ldots,t^{b_r},\epsilon,1,\ldots,1)$. We have
$\alpha\in U\cap S_b$, as desired.
\hfill$\Box$\medskip

We are now in the position to finish the proof of Proposition \ref{cdconnect}. Suppose $S_{[c,d]}$ were not connected. By
Corollary \ref{SEdag} we can write
\begin{equation}
[c,d] = B_1 \coprod B_2\label{eq:decomp}
\end{equation}
and
\begin{equation}
S_{[c,d]}=X_1\coprod X_2,\label{eq:decomp1}
\end{equation}
such that $X_1=S_{B_1}$ and $X_2=S_{B_2}$ are open and closed in
the relative topology of $S_{[c,d]}$. In view of assumption (1) of
Proposition \ref{cdconnect}, (\ref{eq:decomp}) induces a
decomposition of the segment $[c_n,d_n] = B_{1n} \coprod B_{2n}$.

Consider the segment $[c_{n1},d_{n1}]$ of the real line.

\begin{lem}\label{connect} Every point $e_{n1}\in[c_{n1},d_{n1}]$ can
  be covered by a set $V(1)$, open in the topology induced on
  $[c_{n1},d_{n1}]$ by the Euclidean topology of $\R$, having the
  following property. Let $V$ be the subset of $\Gamma^n_+$ defined by
\begin{equation}
V=\{e\in[c,d]\ |\ e_{n1}\in V(1)\}.\label{eq:V}
\end{equation}
Then all of $S_e,e\in V$ are contained in the same set $X_s,\
s=1,2$.
\end{lem}

Lemma \ref{connect} implies Proposition \ref{cdconnect}. Indeed, assume Lemma \ref{connect} and
consider the decomposition (\ref{eq:decomp1}). By Lemma \ref{connect},
(\ref{eq:decomp1}) induces a decomposition $[c_{n1},d_{n1}]=W_1\coprod
W_2$ into two disjoint open sets in the Euclidean topology of $\R$,
which gives the desired contradiction since $[c_{n1},d_{n1}]$ is
connected.

It remains to prove Lemma \ref{connect}.

\noi{\bf Proof of Lemma \ref{connect}:}

\ss{Case 1} : $e_{n1} \notin \sum\limits_{j=1}^{n-1} \Q c_{j1}$. Let
$e_n=(e_{n1},0,\dots,0)$ and $e=(c_1,\dots,c_{n-1},e_n)$. Without loss
of generality, we may assume that $S_e\subset X_1$ (the set $S_e$
consists of a single point). Take a basic open set $U$ of
$\sper\ A$, containing $S_e$ and disjoint from $X_2$. The existence of
$V(1)$ with the desired properties follows immediately from Lemma
\ref{dominant1}.

\ss{Case 2} : $e_{n1} \in \sum\limits_{j=1}^{n-1} \Q c_{j1}$.\smallskip

Let $(e_{n2},\dots,e_{nn})$ be the unique $(n-1)$-tuple of real numbers such that letting
$$
e_n=(e_{n1},e_{n2},\dots,e_{nn}),
$$
we have $e_n\in \sum\limits_{j=1}^{n-1} \Q c_j$.
Let $\pi: A \rightarrow A'$ be a sequence of affine monomial
blowings up with respect to $e$, constructed in Proposition \ref{hir} and Corollary \ref{hir1}, such
that, denoting by $e'$ the transform of $e$ by $\pi$, the real
numbers $e'_{11},\ldots,e'_{n-1,1}$ are $\Q$-linearly independent
and
$$
e'_n=0.
$$
\noi\textbf{Claim.} Assuming $e_n\ne d_n$, we
can choose $\pi$ to be a sequence of blowings up with respect to a
small half-open interval of the form $I=[e,h)\subset[c,d]$ with $h_{n1}>e_{n1}$; in
other words, for any $b\in[e,h)$, and any $\alpha\in S_b$, we have
that $\nu_\alpha$ is non-negative on $A'$.

\noi\textbf{Proof of Claim:} We will go through the procedure of Proposition \ref{hir} step by step and show by induction on
the number of blowings up composing the sequence $\pi$ that there exists a choice of $\pi$ as required in the Claim. Recall
(Remark \ref{Jglobal}) that the algorithm of Proposition \ref{hir} consists of choosing the center of blowing up at each
step. The algorithm works for every choice of coordinate chart at each step, but the choice of the coordinate chart
determines whether the given blowing up is a blowing up with respect to a certain valuation $\nu$.

Let $J$ be the subset of $\{1,\dots,n\}$ prescribed by the algorithm of Proposition \ref{hir} so that the first blowing up
$\pi_1$ appearing in $\pi$ is a blowing up along $(x_J)$.
In the discussion that follows, recall that $c_i=e_i$ for
$i\in\{1,\dots,n-1\}$. Since $\dim_{\mathbb Q}\sum\limits_{i=1}^n\mathbb Qe_{i1}=\dim_{\mathbb Q}\sum\limits_{i=1}^n\mathbb
Qe_i=n-1$, either there exists a unique $j\in J$
such that $e_{j1}<e_{q1}$, $q\in J\setminus\{j\}$, or there exist two indices $j,l\in
J$ such that $e_l=e_j<e_q$, $q\in J\setminus\{l,j\}$. In the first case, we have
\begin{equation}
b_j<b_q,\quad q\in J\setminus\{j\}\ \text{ for every }b=\{b_1,\dots,b_n\}\in[e,h)\label{eq:bj}
\end{equation}
assuming that $h_{n1}-e_{n1}$ is a sufficiently small positive real number. Use $j$ to determine the coordinate chart for
$\pi_1$; then by (\ref{eq:bj}) $\pi_1$ is also a blowing up with respect to every $b\in[e,h)$, as desired.

Next, assume that $e_l=e_j$ are the minimal elements of $e_J$. Take $h\in[e,d]$ with $h_{n1}-e_{n1}>0$ sufficiently small so
that
\begin{equation}
b_l,b_j<b_q,\quad q\in J\setminus\{l,j\}\text{ for every }b=\{b_1,\dots,b_n\}\in[e,h).\label{eq:blbt}
\end{equation}
We have either $b_j\le b_l$ for all $b\in[e,h)$ or $b_j\ge b_l$ for all $b\in[e,h)$ (or both). If $b_j\le b_l$, let $j$
determine the choice of coordinate chart for $\pi_1$; if $b_j\ge b_l$ for $b\in[e,h)$, let the coordinate chart for $\pi_1$
be determined by $l$. In either case, $\pi_1$ is a blowing up with respect to $b$ for all $b\in[e,h)$. Repeating the
above procedure at every step of the algorithm of Proposition \ref{hir}, we can choose $\pi$ as required in the Claim. The
Claim is proved.
\medskip

For any $b\in[e,h)$, we have $S_{b'} \cong S_b$ by
Proposition \ref{valnstays1}. Let $I'$ denote the transform of $I$ under $\pi$, defined in the obvious way.

Since the natural map $\pi^*:\sper\ A' \rightarrow \sper\ A$ is
continuous and preserves rational rank, we will work with
$A',c',d',e',I'$ instead of $A,c,d,e,I$.

The decomposition (\ref{eq:decomp}) induces a decomposition of $I$ and
hence also a decomposition
$$
I'=B'_1\coprod B'_2,
$$
such that $S_{B'_1}$ and $S_{B'_2}$ are open and closed in the induced
topology of $S_{I'}$.

Let $e^\bullet_n=(0,1,0,\dots,0)$,
$e^\bullet=(e'_1,\dots,e'_{n-1},e_n^\bullet)$. Without loss of
generality, assume that $e^\bullet\in B'_1$. Let $U$ be a basic open
set of $\sper\ A$ having non-empty intersection with $S_{e^\bullet}$
and disjoint from $S_{B'_2}$. Now Lemma \ref{dominant2} (with $r=n-1$ and $s=1$) shows that
$S_{e^\bullet} \cap \sur{S_{e'}} \neq \void$ and Lemma
\ref{dominant1} shows that for a sufficiently small interval
$V^\dag(1)=(0,\epsilon)$ of the real line, $\epsilon>0$, if
$b'_{n1}\in V^\dag(1)$ then $S_{b'} \subset U$. Let
$V_+(1)=[0,\epsilon)$ and let $V_+=\left\{e'\in[c',d']\ \left|\
    e'_{n1}\in V_+(1)\right.\right\}$. This proves that
$\bigcup\limits_{b'\in V_+}S_{b'}\subset S_{B'_1}$.

Coming back to the original interval $[c_{n1},d_{n1}]$ (that is,
before performing the blowing up sequence $\pi$), we have shown that all
the $S_b$ for $b\in[e,h)$ lie in the same set $X_s$, $s\in\{1,2\}$,
provided we take $h_{n1}$ sufficiently close to $e_{n1}$ and to the right
of $e_{n1}$. Now, assuming $e_n\ne c_n$, let $J$ be a small half-open
interval of the form $(v,e]$ with $v<e$. Repeating the above reasoning
with $I$ replaced by $J$, we obtain that for $V=(v,h)$,
$\bigcup\limits_{b'\in V}S_{b'}\subset S_{B'_1}$, so the open
interval $V(1)=(v_{n1},h_{n1})$ satisfies the conclusion of the
Lemma.

\noi This completes the proof of Lemma \ref{connect} and Proposition
\ref{cdconnect}.
\hfill $\Box$ \medskip

\noi\textbf{Second proof of Proposition \ref{assumeu=0}:} The main idea of this proof is the following.
Take two points $\alpha$ and $\beta$ in $S_E$. Let
\begin{eqnarray}
a=(a_1,\dots,a_n):&=&(\nu_\alpha(x_1),\dots,\nu_\alpha(x_n))\quad\text{and}\\
b=(b_1,\dots,b_n):&=&(\nu_\beta(x_1),\dots,\nu_\beta(x_n)).
\end{eqnarray}
We join $\alpha$ and $\beta$ by
a ``staircase'' in $\Gamma^n$ where each stair lies entirely in a connected
component of $S_E$. Two examples at the end of the paper show that we have to do this
with some care.\medskip

Let $r\in\{1,\dots,n-1\}$. Consider an element $b\in\Gamma_+^n$ such
that $b_{r+1}=\dots=b_n=0$ and $b_1,\dots,b_r$ are $\Q$-linearly
independent. For $1\le j\le r$, write $b_j=(b_{j1},\dots,b_{jn})$ with
$b_{jq}\in\R$. Consider the $(r\times n)$-matrix $(b_{jq})^{1\le j\le
r}_{1\le q\le n}$.

\begin{lem} There exists $b' =(b'_1,\ldots,b'_n) \in \Gamma^n_+$ having the
following properties:

(1) $b \underset{\circ}{\sim} b'$ (in particular,
$b'_{r+1}=\cdots=b'_n =0$)

(2) writing $b_j'=(b'_{j1},\ldots,b'_{jn})$, we have that the rows of the $(r\times r)$-matrix
  $(b'_{jq})^{1\le j\le r}_{1\le q\le r}$ are $\mathbb Q$-linearly independent.
\end{lem}

\noi \textbf{Proof:} We will gradually replace $b$ by elements of
$\Gamma^n_+$, isomorphic to $b$ in the sense of $\mathcal{OGM}(n)$
until we reach an element of $\Gamma^n_+$ satisfying the conclusion of
the Lemma.

Let $z$ be the number of non-zero columns in the matrix $(b_{jq})^{1\le j\le
r}_{1\le q\le n}$:
$$
z =\#\{q \in \{1,\ldots,n\}\ |\ \exists j,1 \leq j \leq r,\text{ such that }b_{jq}\ne0\}.
$$
Let $\ell$ be the greatest index in $\{1,\ldots,r\}$ such that the column vectors
$$
\left( \begin{array}{c}b_{11}\\ b_{21}\\\vdots \\
b_{r1} \end{array}\right),\ldots,\left(
\begin{array}{c}b_{1\ell}\\ b_{2\ell}\\\vdots\\
b_{r\ell} \end{array}\right)
$$
are $\R$-linearly independent. We proceed by induction on $z-\ell$. First, suppose $z-\ell=0$. Since $\ell=z\le r$, the
$n$-rows $b_1$, \dots, $b_n$ are obtained from the rows of the matrix $(b_{jq})^{1\le j\le r}_{1\le q\le r}$ by
adding $n-r$ zeroes at the end. Since $b_1,\dots,b_r$ are $\Q$-linearly independent by assumption, so are the rows of
$(b_{jq})^{1\le j\le r}_{1\le q\le r}$. Thus in the case $z-\ell=0$ we may take $b'=b$ and there is
nothing more to prove.

Assume $\ell<z$. Permute the columns of $(b_{jq})^{1\le j\le
r}_{1\le q\le n}$ by moving all the zero columns
to the right of the matrix; this does not change the equivalence class of $b$ with respect to the relation $
\underset{\circ}{\sim}$. If this permutation increases the number $\ell$, the proof is finished by induction on $z-\ell$.
Suppose the number $\ell$ stays constant after the above permutation of columns. This means that, after the above
permutation of columnes,
$$
\left(\begin{array}{c}b_{1,\ell+1}\\
b_{2,\ell+1}\\
\vdots\\
b_{r,\ell+1}\end{array}\right)\mbox{ is an $\mathbb R$-linear combination of }
\left(  \begin{array}{c}b_{1j}\\
b_{2j}\\
\vdots\\
b_{rj}
\end{array}\right),\ 1 \leq j \leq\ell.
$$
Replace $b_{j,\ell+1}$ by 0 for $1 \leq j  \leq r$. This operation does not change the equivalence
class of $b$, but decreases the number $z$ of non-zero columns by 1. The result
follows by induction on $z-\ell$.
\hfill$\Box$\medskip

\begin{lem}\label{rr=n} For any connected component $S_{E^\dag}$ of
  $S_E$, there exists $b\in E^\dag$ such that
$$
\rr(b) = n.
$$
\end{lem}

\noi{\bf Proof:} Take $b \in E^\dag$ such that
$\rr(b)=\max\left\{\rr(b^\dag)\ \left|\ b^\dag \in
    E^\dag\right.\right\}$.
Suppose
$$
\rr(b)=r < n.
$$
We will now construct $b^\bullet\in\Gamma_+^n$ such that
\begin{equation}
S_{b^\bullet} \cap \sur{S_b}\neq\void\label{eq:notvoid}
\end{equation}
and
\begin{equation}
\rr(b^\bullet)=r+1.\label{eq:rkr+1}
\end{equation}
This will contradict the maximality of $r$.

Let $A\rightarrow A'$ be a sequence of affine monomial blowings up with respect to $b$,
constructed in \S\ref{Des}, such that, denoting by $b'$ the transform
of $b$ by $\pi$, $b'_1,\ldots,b'_r$ are $\Q$-linearly independent and
$$
b'_{r+1}=\cdots=b'_{n}=0.
$$
Since the natural map $\pi^*:\sper\ A\rightarrow \sper\ A'$ is
continuous and preserves rational rank, we may replace $A$ by $A'$ and
$b$ by $b'$: if the conditions (\ref{eq:notvoid}) and (\ref{eq:rkr+1})
are satisfied in $\sper\ A'$, they will still be satisfied after
applying $\pi^*$ to everything in sight. From now on, we drop the
primes and assume that $b_1,\ldots,b_r$ are $\Q$-linearly independent
and $b_{r+1}=\cdots=b_{n}=0$.

Define
  $b^\bullet=(b_1^\bullet,\dots,b_n^\bullet)$ by
$b^\bullet_i=b_i, i \neq r+1$ and $b^\bullet_{r+1}
=(\underbrace{0,\dots,0}_{r},1,0,\ldots,0)$. By Lemma \ref{dominant2}, (\ref{eq:notvoid}) holds.
Since $b\in E$, we have
\begin{equation}
b^\bullet\in E\label{eq:bbullet}
\end{equation}
because, by the $\Q$-linear independence of the rows of the matrix $(b_{jq})^{1\le j\le
r}_{1\le q\le r}$, the presence of $b^\bullet_{r+1}$ does not affect
the strict inequalities defining $E$: any integer multiple of
$b^\bullet_{r+1}$ is infinitesimal in absolute value compared to any
non trivial $\Z$-linear combination of $b_1,\dots,b_r$. Since $S_b\subset E^\dag$, (\ref{eq:notvoid}) implies that
$S_{b^\bullet}\subset E^\dag$. Since $\rr(b^\bullet)=r+1$, this contradicts the maximality of $r$. This completes
the proof of the Lemma.
\hfill $\Box$ \medskip

For $a=(a_1,\dots,a_n)$, write $a_j=(a_{j1},\dots,a_{jn})$. Take $a
\in E$ such that $a_{11},\ldots,a_{n1}$ are $\Q$-linearly
independent, in particular, $\rr\ a = n$ (such an $a$ exists because
the $n$-tuples $(a_{11},\ldots,a_{n1})$
such that $a_{11},\ldots,a_{n1}$ are $\Q$-linearly independent are
dense in $\R^n$).  \medskip

Now suppose that $S_E$ is not connected and let $S^{(1)}, S^{(2)}$ be
two open and closed sets such that $S_E = S^{(1)} \coprod S^{(2)}$ with
$S^{(1)}$ containing $S_a$. Let $E^{(1)},E^{(2)}$ be the subsets of
$E$ such that $S^{(1)}=S_{E^{(1)}}$ and $S^{(2)}=S_{E^{(2)}}$;
$E^{(1)}$ and $E^{(2)}$ exist by Corollary \ref{SEdag}.

\begin{lem}\label{2nindep}
  There exists $b\in E^{(2)}$ such that
  $a_{11},\ldots,a_{n1},b_{11},\ldots,b_{n1}$ are $\Q$-linearly
  independent.
\end{lem}

\noi{\bf Proof:} According to the preceding lemma, we can find $c \in
E^{(2)}$ such that $\rr(c)=n$. Then $S_c$ consists of a single
point. Let $U$ be a basic open set containing $S_c$ and such that
\begin{equation}
U\cap S^{(1)}=\emptyset.\label{eq:disjoint}
\end{equation}
Let $V\subset\Gamma_+^n$ be a set satisfying the conclusion of Lemma
\ref{dominant1}.
Since the $n$-tuples $(b_{11},\ldots,b_{n1})$ such
that $(a_{11},\ldots,a_{n1},b_{11},\ldots,b_{n1})$ are
$\Q$-linearly independent are dense in $R^n$, there exists
$b\in E\cap V$ such that $(a_{11},\ldots, a_{n1},
b_{11}, \ldots, b_{n1})$ are $\Q$-linearly
independent. By Lemma \ref{dominant1}, $S_b\subset U$. By
(\ref{eq:disjoint}), $b\in E^{(2)}$, as desired.
\hfill $\Box$ \medskip

\noi Let $a \in E^{(1)}$, $b \in E^{(2)}$ be such that
$a_{11},\ldots,a_{n1},b_{11},\ldots,b_{n1}$ are $\Q$-linearly independent.
\begin{rek}  If $\lambda,\mu \in \Q,\ \lambda \neq \mu$, then
$\lambda a_{11} + (1-\lambda) b_{11}$, $\ldots,\lambda a_{n1} +
(1-\lambda) b_{n1}$, $\mu a_{11} + (1-\mu) b_{11}$, \ldots, $\mu a_{n1}+
(1-\mu) b_{n1}$ are $\Q$-linearly independent.
\end{rek}

In the sequel, let $N$ be a large natural number and $\dsp{\lambda =
  \frac{i}{N},\ \mu= \frac{i+1}{N}}$.
\begin{lem} For $N \in \N$ sufficiently large we have
$$
\left( \begin{array}{c} \left(\frac{i}{N}
    a_{11} + \left(1-\frac{i}{N}\right)b_{11},0,\ldots,0\right) \\ \vdots \\
    \left(\frac{i}{N} a_{j-1,1} +
    \left(1-\frac{i}{N}\right)b_{j-1,1},0,\ldots,0\right) \\
    \left(\frac{i+1}{N}a_{j1} + \left(1-\frac{i+1}{N}\right)b_{j1},0,\ldots,0\right) \\
    \vdots \\
    \left(\frac{i+1}{N} a_{n1} +
    \left(1-\frac{i+1}{N}\right)b_{n1},0,\ldots,0\right) \end{array} \right) \in E.
$$
for all $i \in \{0,\ldots,N-1\}$ and all $j\in \{1,\ldots,n\}$.
\end{lem}
\noi{\bf Proof:} Since $a_{11},\ldots,a_{n1}$ are $\Q$-linearly
independent, saying that $ \sum\limits_{i=1}^n \omega_{il}a_i > 0, \ l
\in \{1,\ldots,q\}$ is equivalent to saying that
$$
\sum_{i=1}^n \omega_{il}a_{i1} > 0, \ l \in \{1,\ldots,q\}.
$$
Similarly, we have
$$
\sum_{i=1}^n \omega_{il}b_{i1} > 0, \ l \in \{1,\ldots,q\}.
$$
The set of $n$-tuples $(c_{11},\dots,c_{n1})\in\R^n$ such that
$\sum\limits_{i=1}^n \omega_{il}c_{i1} > 0, \ l \in \{1,\ldots,q\}$, is open and convex and contains both
$(a_{11},\dots,a_{n1})$ and $(b_{11},\dots,b_{n1})$.
Then for $N$ sufficiently large we have
$$
\sum\limits_{i=1}^{j-1}\omega_{il}\left(\frac
iNa_{i1}+\left(1-\frac iN\right)b_{i1}\right)+\sum\limits_{i=j}^n\omega_{il}\left(\frac{i+1}Na_{i1}+
\left(1-\frac{i+1}N\right)b_{i1}\right)=
$$
$$
\sum\limits_{i=1}^n\omega_{il}\left(\frac iNa_{i1}+\left(1-\frac
iN\right)b_{i1}\right)+\sum\limits_{i=j}^n\omega_{il}\left(\frac1Na_{i1}-\frac1Nb_{i1}\right)>0,
$$
$l\in \{1,\ldots,q\}$. This completes the proof of the Lemma.
\hfill$\Box$\medskip

For each $i \in\{0,\ldots,N-1\}$ and $j \in \{1,\ldots,n\}$, consider the pair
\begin{eqnarray}
a(i,j-1):=&\left(
  \begin{array}{c} \left(\frac{i}{N}a_{11} +
    \left(1-\frac{i}{N}\right)b_{11},0,\ldots,0\right) \\ \vdots \\
    \left(\frac{i}{N}a_{j-1,1} + \left(1-\frac{i}{N}\right)
    b_{j-1,1},0,\ldots,0\right) \\
    \left(\frac{i+1}{N}a_{j1} + \left(1-\frac{i+1}{N}\right)b_{j1},0,\ldots,0\right) \\
    \vdots \\
    \left(\frac{i+1}{N}a_{n1} + \left(1-
    \frac{i+1}{N}\right)b_{n1},0,\ldots,0\right) \end{array}
\right),\\
a(i,j):=&
\ \ \left( \begin{array}{c} \left(\frac{i}{N}a_{11} +
    \left(1-\frac{i}{N}\right)b_{11},0,\ldots,0\right) \\
    \vdots \\
    \left(\frac{i}{N}a_{j,1} + \left(1-\frac{i}{N}\right)
    b_{j,1},0,\ldots,0\right) \\
    \left(\frac{i+1}{N}a_{j+1,1} + \left(1-\frac{i+1}{N}\right)b_{j+1,1},0,\ldots,0\right)
    \\ \vdots \\
    \left(\frac{i+1}{N}a_{n1} +\left(1-\frac{i+1}{N}\right)b_{n1},0,\ldots,0\right) \end{array} \right).
\end{eqnarray}
of elements of $\Gamma_+^n$.

We can now finish proving the connectedness of $S_E$. Take $i\in\{0,\ldots,N-1\}$ and
$j\in\{1,\ldots,n\}$. Up to renumbering the components
$c_1,\dots,c_n$, the pair of points $c=a(i,j-1)$ and $d=a(i,j)$
satisfy the hypotheses of Proposition \ref{cdconnect}. Then
Proposition \ref{cdconnect} says that $S_{a(i,j-1)}$ and
$S_{a(i,j)}$ are contained in the same connected component of
$S_E$, so $a(i,j-1)$ and $a(i,j)$ belong to the same set $E^{(s)},\
s=1,2$. Put
$$
b_{trunc}:=\left(
\begin{array}{c} (b_{11},0,\ldots,0) \\ \vdots \\
(b_{n1},0,\ldots,0) \end{array}
\right).
$$
Since $b\in E^{(2)}$ and $S_b=S_{b_{trunc}}$ by linear
independence of $b_{11},\dots,b_{n1}$, we must have $b_{trunc}\in E^{(2)}$.

Since $b_{trunc}=a(0,n)$ and $a(i,j-1)$ and $a(i,j)$ belong to the
same set $E^{(s)}$ for all $i,j$, we have $a(i,j)\in E^{(2)}$ for all
$i,j$ by induction on $(i,n-j)$ in the lexicographical ordering (note
that $a(i,0)=a(i+1,n)$). Then
$$
a(N-1,0)=a_{trunc}:=\left(
  \begin{array}{c} (a_{11},0,\ldots,0) \\ \vdots \\
     (a_{n1},0,\ldots,0) \end{array}
\right)\in E^{(2)}
$$
so, finally, $a\in E^{(2)}$ which contradicts the fact that $a$ was
chosen to lie in $E^{(1)}$. This completes the second proof of Proposition \ref{assumeu=0} (the connectedness of $S_E$) and
with it of Theorem \ref{con1}.\medskip

Next, we explain how to use Theorem \ref{con1} to give another proof of Theorem \ref{conmonomial}. Let the notation
be as in the proof of Theorem \ref{conmonomial}. Let $C$ be the set of all $\delta\in\sper\ A$ such that
the inequalities (\ref{eq:ineqalpha2}) are satisfied with $\alpha$ replaced by $\delta$ and $x_q(\delta)$ has the same sign
as $x_q(\alpha)$ for all $q$ such that $x_q\notin<\alpha,\beta>$. The set $C$ contains both $\alpha$ and $\beta$ and none
of the polynomials $g_i$ change sing on $C$. Thus to complete the proof of Theorem \ref{conmonomial}, it remains to prove
\begin{cor} (of Theorem \ref{con1}). The set $C$ is connected.
\end{cor}
\textbf{Proof:} For an index $q\le l$, replacing $x_q$ by $-x_q$ does not change the problem. Thus we may assume that the
sign conditions appearing in the definition of $C$ are all of the form $x_q>0$ for $q\le l$. We proceed by induction on
$n-l$. Il $l=n$ then the Corollary is a special case of Theorem \ref{con1}. Assume that $l<n$ and that the result is known
for all the smaller values of $n-l$. Let $C_+=\{\delta\in C\ |\ x_n(\delta)>0\}$, $C_-=\{\delta\in C\ |\ x_n(\delta)<0\}$
and $C_0=\{\delta\in C\ |\ x_n(\delta)=0\}$. By the induction assumption, each of $C_+$, $C_-$ and $C_0$ is connected.
Take any point $\delta\in C_0$ and a basic open set $U$ of $\sper\ A$, containing $\delta$. It is easy to show, by an
argument similar to that of Lemma \ref{dominant2}, that $U\cap C_+\ne\emptyset$ and $U\cap C_-\ne\emptyset$. Hence
$C_0\subset\bar C_+$ and $C_0\subset\bar C_-$. This proves that $C=C_+\cup C_0\cup C_-$ is connected, as desired.
\hfill$\Box$\medskip

The next two examples show that the intuition that ``convex sets in
$\Gamma^n_+$ give rise to connected sets in $\sper\ A$'' is not
completely accurate.
\medskip

\noi \textbf{Example 1:} Let $n=2$. Let
$c=\left(\begin{array}{cc}0&0\\0&0\end{array}\right),\ d=
\left(\begin{array}{cc}0&0\\0&1\end{array}\right) \in\ \Gamma^2$. Let $S_{[c,d]}$ be the subset of $\sper\ A$, defined in
(\ref{eq:[c,d]}) and (\ref{eq:Scd}). Let $S_{]c,d]}=S_{[c,d]}\setminus S_c$. Then $S_{[c,d]}$ is not connected:
we have $S_{[c,d]}=S_c\coprod S_{]c,d]}=S_c \coprod S_d$ and both $S_c$ and
$S_d$ are open and closed. In particular, Proposition
\ref{cdconnect} is false without the assumption that
$c_{11},\ldots,c_{n1},d_{n1}$ are linearly independent. Another way of
interpreting this example is that intervals consisting only of rank 1
valuations do not give rise to connected sets in $\sper\ A$.
\medskip

\noi\textbf{Example 2:} In this example, higher rank valuations
appear, but, again, $c_{11},\ldots,c_{n1},d_{n1}$ are $\Q$-linearly
dependent. Let $n=3$ and take $c_1=d_1=e_1=(1,0,0)$,
$c_2=d_2=e_2=(0,1,0)$, $c_3 =(0,0,0)$, $d_3=(1,0,0)$,
$e_3=(e_{31},e_{32},e_{33}),\ c_3 \leq e_3 \leq d_3$. Then
$c_{11},c_{21},c_{31},d_{31}$ are $\Q$-linearly dependent and the
segment $[c_{31},d_{31}]=[0,1]$ is a disjoint union of $V^{(1)}=
\{e_{31}=0\}$ and $V^{(2)}=\{0<e_{31}\le1\}$. The decomposition
$[c_{31},d_{31}]=V^{(1)}\coprod V^{(2)}$ induces a decomposition
$S_{[c,d]}=S_{V^{(1)}}\coprod S_{V^{(2)}}$ into two disjoint
sets. We claim that $S_{V^{(1)}}$ and $S_{V^{(2)}}$ are both open and closed in the induced topology of $S_{[c,d]}$. To see
this, it is sufficient to cover $S_{V^{(1)}}$ by basic open sets of $\sper\ A$, disjoint from $S_{V^{(2)}}$, and vice versa.

Let $\N_+$ denote the set of strictly positive integers.
For each $q\in\N_+$, let $U_q$ denote the basic open set of $\sper\ A$, defined by the equation $x_3-x_2^q>0$. For a
point $\delta\in S_{V^{(2)}}$ we have $q\nu_\delta(x_2)=(0,q,0)<(e_{31},e_{32},e_{33})=\nu_\delta(x_3)$ (see the definition
of the lexicographical ordering at the end of section \S\ref{Cu}), so $x_2^q>_\delta x_3$ and $\delta\notin U_q$. On the
other hand, if $\delta\in S_{V^{(1)}}$ then $\nu_\delta(x_3)$ has the form $(0,e_{32},e_{33})$. Taking $q>e_{32}$, we have
$q\nu_\delta(x_2)=(0,q,0)>(0,e_{32},e_{33})=\nu_\delta(x_3)$, so $x_2^q<x_3$ and $\delta\in U_q$. Thus
$S_{V^{(1)}}\subset\bigcup\limits_{q\in\N_+}U_q$ and $U_q\cap S_{V^{(2)}}=\emptyset$ for all $q\in\N_+$. This proves that
$S_{V^{(1)}}$ is relatively open in the induced topology.

Next, let $W_q$ denote the basic open set of $\sper\ A$, defined by the equation $x_1 - x_3^q>0$, $q \in \N_+$. For a point
$\delta\in S_{V^{(1)}}$ we have $q\nu_\delta(x_3)=(0,qe_{32},qe_{33})<(1,0,0)=\nu_\delta(x_1)$, so $x_3^q>x_1$ and
$\delta\notin W_q$. On the other hand, if $\delta\in S_{V^{(2)}}$ then
$\nu_\delta(x_3)$ has the form $(e_{31},e_{32},e_{33})$ with $e_{31}>0$. Taking $q>\frac1{e_{31}}$, we have
$q\nu_\delta(x_3)=(qe_{31},qe_{32},qe_{33})>(1,0,0)=\nu_\delta(x_1)$, so $x_3^q<x_1$ and $\delta\in W_q$. Thus
$S_{V^{(2)}}\subset\bigcup\limits_{q\in\N_+}W_q$ and $W_q\cap S_{V^{(1)}}=\emptyset$ for all $q\in\Q_+$. This proves that
$S_{V^{(2)}}$ is relatively open in the induced topology, as desired.

\end{document}

%% file: connfig.pstex_t
\begin{picture}(0,0)%
\includegraphics{connfig.pstex}%
\end{picture}%
\setlength{\unitlength}{2735sp}%
\begingroup\makeatletter\ifx\SetFigFont\undefined%
\gdef\SetFigFont#1#2#3#4#5{%
  \reset@font\fontsize{#1}{#2pt}%
  \fontfamily{#3}\fontseries{#4}\fontshape{#5}%
  \selectfont}%
\fi\endgroup%
\begin{picture}(5481,2628)(3463,-3023)
\put(4411,-1636){\makebox(0,0)[lb]{\smash{{\SetFigFont{8}{9.6}{\rmdefault}{\mddefault}{\updefault}{\color[rgb]{0,0,0}$\alpha$}%
}}}}
\put(8416,-1411){\makebox(0,0)[lb]{\smash{{\SetFigFont{8}{9.6}{\rmdefault}{\mddefault}{\updefault}{\color[rgb]{0,0,0}$\beta$}%
}}}}
\put(7471,-2446){\makebox(0,0)[lb]{\smash{{\SetFigFont{8}{9.6}{\rmdefault}{\mddefault}{\updefault}{\color[rgb]{0,0,0}$f_{i_{s-1}}$}%
}}}}
\put(8191,-2086){\makebox(0,0)[lb]{\smash{{\SetFigFont{8}{9.6}{\rmdefault}{\mddefault}{\updefault}{\color[rgb]{0,0,0}$f_{i_s}$}%
}}}}
\put(5311,-2671){\makebox(0,0)[lb]{\smash{{\SetFigFont{8}{9.6}{\rmdefault}{\mddefault}{\updefault}{\color[rgb]{0,0,0}$f_{i_2}$}%
}}}}
\put(4456,-2446){\makebox(0,0)[lb]{\smash{{\SetFigFont{8}{9.6}{\rmdefault}{\mddefault}{\updefault}{\color[rgb]{0,0,0}$f_{i_1}$}%
}}}}
\end{picture}%